\documentclass[12pt]{article}

\usepackage{amsmath}
\usepackage{amsfonts}
\usepackage{amssymb}
\usepackage{amsthm}
\usepackage{graphicx}

\usepackage{cite}
\usepackage{caption}
\usepackage{url}
\usepackage{color}\usepackage{placeins}
\usepackage{svg}
\usepackage{algorithm}
\usepackage{algpseudocode}
\usepackage{mathtools}
\usepackage{subcaption}
\usepackage{multirow}
\usepackage{bbm}
\usepackage{epstopdf}
\usepackage{graphicx}
\usepackage{bm}
\usepackage{enumerate}
\usepackage{booktabs}

\theoremstyle{definition}

\newtheorem{remark}{Remark}
\newcommand{\norm}[1]{\left\lVert#1\right\rVert}

\algdef{SE}{Begin}{End}{\textbf{begin}}{\textbf{end}}

\title{General-domain FC-based shock-dynamics solver II:\\ Non-smooth
  domains, accuracy and parallel performance}

\author{Daniel V. Leibovici\footnote{NVIDIA Research, Santa Clara, CA, 95051, USA} \and Oscar P. Bruno\footnote{Computing and Mathematical Sciences, Caltech,
    Pasadena, CA 91125, USA}}
    
\begin{document}
\date{}
\maketitle


\begin{abstract}
  This contribution, Part II of a two-part series, extends the
general-domain FC-SDNN (Fourier Continuation Shock-detecting Neural
Network) introduced in Part~I to enable treatment of non-smooth
domains, it introduces a parallel implementation of the scheme with
high-quality weak and strong scalability properties, and it
illustrates the overall methodology with a variety of tests for the 2D
Euler equations---including supersonic and hypersonic flows and shocks
past obstacles with corners. The results produced by the new methods
are compared to previous theoretical and experimental results, and the
high parallel scalability of the algorithm is demonstrated in both
weak and strong scalability tests.  Thanks to its use of a localized
yet smooth artificial viscosity term---whose support is confined to
regions near flow discontinuities identified by an artificial neural
network---the algorithm maintains minimal numerical dissipation away
from discontinuities. Overall, the method delivers accurate, sharp
resolution of shocks and contact discontinuities, while producing
smooth numerical solutions in regular flow regions—as evidenced by the
near-complete absence of spurious oscillations in level set contours,
even under strong shocks and high-speed flow conditions.

\end{abstract}

\vspace{0.5 cm}
\maketitle
\newpage
\section{\label{sec:introduction}Introduction}

This article is the Part~II in a two paper sequence devoted to the
development of a general-domain {\em FC-SDNN} (Fourier Continuation
Shock-detecting Neural Network) spectral scheme for the numerical
solution of nonlinear conservation laws under arbitrary boundary
conditions in {\em general domains}.

As outlined in the introduction of the Part~I
contribution~\cite{bruno_leibo_partI}, the FC-SDNN algorithm
(i)~Discretizes the flow variables by means of an overlapping-patch
domain decomposition approach and patchwise spatial Fourier
continuation-based spectral representation of non-periodic flow
variables on quasi-uniform meshes (which result in limited dispersion
and mild CFL restrictions); (ii)~Detects shocks and other flow
discontinuities on the basis of a Gibbs-ringing neural-network
smoothness classifier; (iii)~Effectively eliminates spurious
oscillations near shocks and other discontinuities while preserving
sharp resolution of flow features and maintaining high accuracy in
smooth flow regions, thanks to its carefully designed smooth and
localized artificial viscosity assignments; and, (iv)~Does not require
the use of positivity preserving limiters in view of a setting that,
consistent with Prandtl's boundary-layer theory, and unlike other
existing solvers, natively incorporates boundary-layer structures near
physical boundaries along with inviscid flow away from such
boundaries.  The present Part-II, in turn, (v)~Provides the
algorithmic extensions needed for the treatment of domains with
corners; (vi)~Introduces a parallel implementation of high parallel
efficiency and scalability; and (vii)~It validates and benchmarks the
method against a variety of existing experimental and theoretical
results.

The FC-based solution of PDEs in domains with corners requires a new
strategy, as the previously introduced FC discretizations are
inherently suited to spatial domains formed by overlapping regions,
each representing the image of a smooth parametrization over a square
domain---a discretization approach that cannot be directly extended to
include corner points. Indeed, a seemingly straightforward
approach---using smooth patches on each side of the corner that extend
smoothly beyond the corner to ensure overlap---leads to
instabilities. This difficulty arises because such smoothly continued
patches lead to decompositions near corners that fail to satisfy a
minimum overlap condition, which, as demonstrated by extensive
experiments conducted in the course of this work, is essential for
stability. Specifically, the minimum overlap condition requires each
point on a patch boundary to lie sufficiently deep within the interior
of another patch; see Part-I, Section 3.2.2 for smooth domains and
Section~\ref{subsec:overlap} in this paper for non-smooth domains. In
order to avoid this difficulty this paper develops the concept of
``corner patches'' and, specifically, corner patches of two types: one
for convex corners and one for concave corners. The use of such
specialized corner-patch subdomains ensures the necessary stability
while effectively addressing the geometric constraints posed by
corners in FC-based gas-dynamics solutions. To the best of our
knowledge, the spectral solvers presented in this two-part sequence
are the first to handle shocks in general domains: previous spectral
solvers, such as those in~\cite{bruno2022fc,
  guermond2011entropy,kornelus2018flux}, or hybrid physics-ML solvers
in which corrections to a spectral solver are learned using a neural
network~\cite{dresdner2022learning}, were demonstrated only for flows
in rectangular domains without obstacles.

The overlapping-patch decomposition approach just described not only
provides geometric and algorithmic flexibility but also seamlessly
integrates with parallelization strategies for the efficient solution
of shock-dynamics problems. For reference we note that various methods
have previously been used to achieve parallelization of shock-dynamics
solvers, including multi-domain decompositions based on use of ghost
points for WENO schemes~\cite{chao2009massively} as well as
parallelization strategies appropriate for finite-element
discretizations for the Discontinuous Galerkin
methods~\cite{baggag1999parallel, sonntag2017efficient} and the
Entropy viscosity approach~\cite{ryujin-2021-1, ryujin-2021-3}. 


The proposed spectral, multi-patch strategy is validated through
various tests, including a comparison with
the single-patch case from~\cite{bruno2022fc}
(Section~\ref{validation}) and an analysis of the scheme’s energy
conservation properties (Section~\ref{engy_cons}), as well as
comparisons with exact solutions for supersonic and hypersonic flow
past a triangular wedge (Section~\ref{subsubsec:wedge}), test problems
considered in previous contributions (Sections~\ref{subsubsec:wedge},
\ref{subsubsec:flow_prism}, and~\ref{Matrixcylinders}), and
illustrations in the context of existing experimental data
(Section~\ref{subsubsec:shock_prism}).

This paper is organized as follows. After a brief overview in
Section~\ref{Preliminaries} of the multi-patch FC-SDNN Euler solver
introduced in~Part~I, Sections~\ref{sec:Multidomain}
and~\ref{sec:viscosity} describe the FC-SDNN overlapping patch
strategy and associated multi-patch artificial viscosity assignment,
respectively, in the context of non-smooth domains considered in this
paper. An overview of the overall multi-patch FC-SDNN algorithm for
non-smooth domains is then presented in Section~\ref{main-elms}, in
preparation for the introduction of the proposed parallel multi-domain
strategy and the overall FC-SDNN MPI implementation in
Section~\ref{sec:Parallel}. A number of test cases presented in
Section~\ref{sec:numerical results} demonstrate the favorable strong
and weak scalability properties of the parallel FC-SDNN scheme,
validate the multi-patch strategy, and illustrate the method's
performance for supersonic and hypersonic flows and shocks past
non-smooth obstacles. Section~\ref{sec:conclusion} presents a few
concluding remarks. 
\begin{remark}\label{referncs}
  Throughout this paper, references to equations and sections in the
  Part~I contribution are indicated by preceding the equation or
  section number by the roman numeral 'I' so that, e.g., equation (5)
  and Section 3.2 in Part~I are referred to as equation I--(5) and
  Section~I--3.2 respectively.
\end{remark}

\section{\label{Preliminaries} Preliminaries and 
  Part~I background}

Following Part~I, this paper proposes novel Fourier spectral
methodologies for the numerical solution of the two-dimensional Euler
equations in general 2D non-periodic domains $\Omega$ and with general
boundary conditions; unlike Part~I, the algorithms presented in this
paper allow for consideration of domains $\Omega$ that incorporate
obstacles containing concave and/or convex corners. In detail, we
consider the 2D Euler equations
\begin{equation} \label{eq: euler 2d equation}
\dfrac{\partial}{\partial t}
\begin{pmatrix}
\ \rho \\[\jot]
\ \rho u\\[\jot]
\ \rho v\\[\jot]
\ E
\end{pmatrix} + \frac{\partial}{\partial x}
\begin{pmatrix}
\ \rho u\\[\jot]
\ \rho u^2 + p\\[\jot]
\ \rho u v\\[\jot]
\ u (E + p)
\end{pmatrix}
+ \frac{\partial}{\partial y}
\begin{pmatrix}
\ \rho v\\[\jot]
\ \rho u v\\[\jot]
\ \rho v^2 + p\\[\jot]
\ v (E + p)
\end{pmatrix} = 0
\end{equation}
on the domain $\Omega\subset \mathbb{R}^2 $, where $\rho$,
\textit{\textbf{u}}, $p$ denote the density, velocity
vector, and pressure respectively; and the total energy $E$ is given by
\begin{equation}
  E = \frac{p}{\gamma - 1} + \frac{1}{2}\rho \lvert \textit{\textbf{u}} \rvert ^2,
\end{equation}
where $\gamma$ denotes the heat capacity ratio; the ideal diatomic-gas
heat-capacity value $\gamma = 1.4$ is used in all of the test cases
considered in this paper. As in Part~I and with reference to
Remark~\ref{referncs} above, in order to incorporate adiabatic
boundary conditions, the well-known expression~I--(4) that expresses
the energy $E$ in terms of the temperature $\theta$ together with
$\rho$, $\textit{\textbf{u}}$ and $\gamma$ is utilized in this
paper. Clearly equation~\eqref{eq: euler 2d equation} may be expressed
in the form
\begin{equation} \label{eq: nonlinear-eqn} \dfrac{\partial}{\partial
    t} \textbf{e}(\textbf{x}, t) + \nabla \cdot \big(
  f(\textbf{e}(\textbf{x}, t))\big) = 0,\quad \textbf{e} : \Omega \times [0, T] \rightarrow \mathbb{R}^4,
\end{equation}
for a certain function
$f : \mathbb{R}^4 \rightarrow \mathbb{R}^4 \times \mathbb{R}^2$ and
vector $\textbf{e}= \textbf{e}(\textbf{x}, t)$ displayed in
equation~I--(5).

The FC-SDNN algorithm relies on four main elements, namely, (i)~The
FC-Gram algorithm for spectral representation of non-periodic
functions; (ii)~A multi-patch domain-decomposition strategy; (iii)~A
Shock-Detecting Neural Network (SDNN) for assignment of artificial
viscosity; and (iv)~An accurate time-marching scheme that additionally
does not introduce oscillations in presence of shocks and other flow
discontinuities; in what follows we briefly review these algorithmic
elements which are otherwise discussed at length in Part~I and
references therein.

The FC-Gram algorithm mentioned in point~(i) constructs an accurate
Fourier approximation of a given, generally non-periodic, function $F$
defined on a given one-dimensional interval, which, without loss of
generality, is assumed in this section to equal the unit interval, so
that $F:[0,1]\to \mathbb{R}$. Starting from given function values
$F_j =F(x_j)$ at $N$ points $x_j=jh\in [0,1]$ ($h=1/(N-1)$,
$j=0,\dots,N-1$), the FC-Gram algorithm produces a function $F^c(x)$
given by a Fourier expansion containing a number $M\sim N/2$ of
Fourier modes, which is defined (and periodic) in an interval
$[0,\beta]$ that strictly contains $[0,1]$, and which closely
approximates $F(x)$ for $x\in [0,1]$. Details concerning the
construction of the Fourier expansion $F^c(x)$, as well as the
rationale underlying the use of the FC-Gram algorithm of degree $d=2$,
which is a common feature of both the Part~I and the present Part~II
contributions, can be found in Section~I--2.2; here we merely mention
that, as indicated e.g in Part~I and detailed
in~\cite[Sec. 4.1.2]{bruno2022fc} and references therein, use of the
FC-Gram algorithm leads to PDE solvers of low dissipation and
dispersion---as it behooves a spectrally-based discretization
methodology.

Like Part~I and per point~(ii) above, the shock-dynamics solver
presented in this paper relies on the use of decomposition of the overall
computational domain as a union of overlapping patches. The character
of the decomposition is analogous to the one used in Part~I, as it
includes patches of various types covering regions away from
boundaries and near boundaries. Unlike the decompositions considered
in Part~I, however, all of which are obtained as images of
parametrizations defined in the unit-square parameter space, the
present paper introduces, in addition, certain patch types around
corners of any existing obstacle, for which the parameter spaces are
given either by the unit square or by a prototypical L-shaped
domain. A description of all patches used, including those considered
in Part~I, as well as the new $\mathcal{C}_1$- and
$\mathcal{C}_2$-type corner patches introduced in this paper for
convex and concave curvilinear corners respectively, and associated
discretization and inter-patch solution communication procedures, are
presented in Section~\ref{sec:Multidomain}.

Following upon~\cite{schwander2021controlling,bruno2022fc} as well as
Part~I, the algorithm introduced in this paper relies on the use of a
shock detecting neural network and associated viscosity assignments to
control spurious oscillations emanating from shocks and other flow
discontinuities. The multi-patch artificial viscosity assignment
method considered here---which is based on the use of subpatch-wise
produced viscosity values (obtained, at each space-time point in each
patch from an expression involving the meshsize and the maximum wave
speed at that point) that are subsequently combined by means of an
inter-patch windowed-viscosity propagation algorithm across patches
and subpatches---is actually identical to the one presented in Section
I--5 and it will therefore not be described in this paper at any
length; all relevant details can be found in the Part~I section just
mentioned. Similarly, the time-stepping algorithm coincides with the
one introduced in Section~4 in Part~I, to which we refer for
details. In brief, the algorithm time-steps the solution by means of
an SSPRK-4 scheme~\cite{gottlieb2005high} on each subpatch, enforcing
boundary conditions at each stage of the scheme, with subpatch-wise
solution filtering at each time-step (to control error growth in
unresolved high frequency modes, as befits traditional spectral
methods), together with a localized discontinuity smearing strategy
for initial data which is applied before the first
time-step. Importantly, as detailed in Remarks~3 and~4 in Part~I, the
algorithm enforces full viscous boundary conditions at obstacle
boundaries in conjunction with adiabatic boundary conditions for the
temperature---which, since the temperature is not one of the
conservative variables used, is enforced indirectly by expressing the
internal energy $E$ in terms of temperature, density and velocity.

\begin{figure}[H]
\centering
\includegraphics[width=0.5\linewidth,]{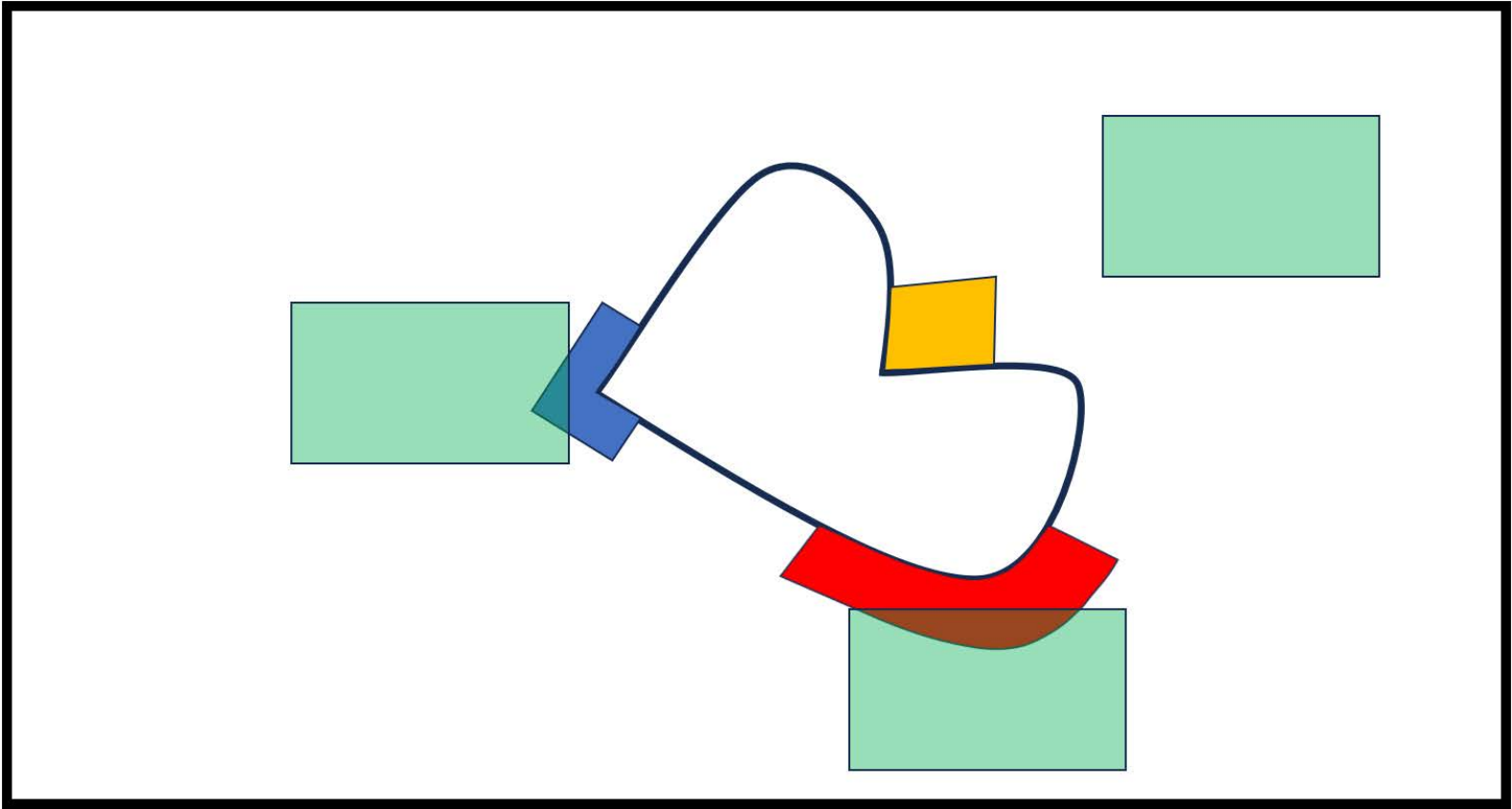}
\caption{Patch types used in the overlapping-patch decomposition of a
  general domain $\Omega$, including (a)~A corner patch of type
  $\mathcal{C}_2$ (orange); (b)~A smooth-boundary patch of type
  $\mathcal{S}$ (red); (c)~Interior patch of type $\mathcal{I}$
  (green); and (d)~A corner patch of type $\mathcal{C}_1$ (blue). The
  $\mathcal{Q}$ parameter space is used for~(a) through~(c) and the
  $\mathcal{L}$ parameter space is used for~(d).}
    \label{fig:mappings}
\end{figure}

\section{\label{sec:Multidomain} Overlapping-patch/subpatch geometry description and discretization}
\subsection{\label{subsec:decomposition} Domain decomposition: patches and their parametrizations}

Following~\cite{chesshire1990composite,albin2011spectral,amlani2016fc,bruno2019higher},
in order to obtain an FC-based discretization of equation~\eqref{eq:
  euler 2d equation} for general 2D domains $\Omega$ which may feature
curved and/or non-smooth boundaries $\Gamma = \partial \Omega$, the
proposed method relies on the decomposition of the connected open set
$\Omega$ as a union of a number
$P = P_{\mathcal{S}} + P_{\mathcal{C}_1}+ P_{\mathcal{C}_2}+
P_{\mathcal{I}}$ of overlapping patches (see Remark~\ref{rem_bdries}),
each one of which is an open set, including, a number
$P_{\mathcal{S}}$ of $\mathcal{S}$-type patches
$\Omega^{\mathcal{S}}_{p}$ (smooth-boundary patches); a number
$P_{\mathcal{C}_1}$ of $\mathcal{C}_1$-type patches
$ \Omega^{\mathcal{C}_1}_{p}$; a number $P_{\mathcal{C}_2}$ of
$\mathcal{C}_2$-type patches $ \Omega^{\mathcal{C}_2}_{p}$; and a
number $P_{\mathcal{I}}$ of $\mathcal{I}$-type patches
$\Omega^{\mathcal{I}}_{p}$ (interior patches):
\begin{equation}\label{decomp}
  \Omega = \left(\bigcup_{p=1}^{P_{\mathcal{S}}} \Omega^{\mathcal{S}}_{p}\right ) \cup \left(\bigcup_{p=1}^{P_{\mathcal{C}_1}} \Omega^{\mathcal{C}_1}_{p}\right )  \cup \left(\bigcup_{p=1}^{P_{\mathcal{C}_2}} \Omega^{\mathcal{C}_2}_{p}\right ) \cup \left(\bigcup_{p=1}^{P_{\mathcal{I}}} \Omega^{\mathcal{I}}_{p}\right ).
\end{equation}
Here $\mathcal{C}_1$-type patches (resp. $\mathcal{C}_2$-type patches)
are defined to contain a neighborhood within $\Omega$ of boundary
corner points of ``Type $\mathcal{C}_1$'' (also called $\mathcal{C}_1$
corner points in what follows), namely corner points formed by two
smooth arcs meeting at $C$ with interior angles $\theta > 180^\circ$
(resp. corner points of ``Type $\mathcal{C}_2$'', or $\mathcal{C}_2$
corner points, namely corner points with interior angles $\theta$
satisfying $0 <\theta < 180^\circ$). $\mathcal{S}$-type patches and
$\mathcal{I}$-type patches, on the other hand are used to cover areas
along smooth boundary regions and regions away from boundaries,
respectively.

\begin{remark}\label{CeqC2}
  Note that the $\mathcal{S}$- and $\mathcal{I}$ patch-types in the
  domain decomposition~\eqref{decomp} are identical to the homonymous
  patch types introduced in Part I. The $\mathcal{C}_2$ patch-type
  included in the decomposition~\eqref{decomp}, in turn, generalizes
  the $90^\circ$ interior-angle $\mathcal{C}$ patch-type introduced in
  Part~I: as illustrated in the middle top panel in
  Figure~\ref{fig:corner_mapping}, the $\mathcal{C}_2$ patch-type
  corresponds to patches around a corner point, with possibly curved
  sides, and with interior angle less than $180^\circ$. The present
  contribution, finally, introduces the new $\mathcal{C}_1$-type
  patch, illustrated in the leftmost top panel in
  Figure~\ref{fig:corner_mapping}, whose interior angles are larger
  than $180^\circ$.

\end{remark}

\begin{remark}\label{rem_bdries}
  As detailed in the following sections, each patch
  $\Omega^{\mathcal{R}}_{p}$ ($\mathcal{R}=\mathcal{S}$,
  $\mathcal{I}$, $\mathcal{C}_1$, $\mathcal{C}_2$;
  $p=1,\dots,P_{\mathcal{R}}$) equals the interior of the image under a
  smooth parametrization of a closed polygonal parameter domain
  $\mathcal{P}$. A requirement imposed on the patches
  $\Omega^{\mathcal{R}}_{p}$, which is satisfied in all cases by the
  patches constructed in the following sections, is that the image of
  each one of the sides of the polygonal parameter domain is either fully
  contained within $\Gamma$ or it contains at most one point of
  intersection with $\Gamma$ (the latter case occurring as the image
  of one of the parameter polygon sides meets $\Gamma$ transversally.)
\end{remark}

As detailed in what follows, $\mathcal{S}$-, $\mathcal{I}$- and
$\mathcal{C}_2$-type patches (resp.$\mathcal{C}_1$-type patches) are
discretized on the basis of smooth invertible mappings, each one of
which maps the closed canonical parameter space onto the closure
$\overline{\Omega^\mathcal{R}_{p}}$ of a single open patch
$\Omega^\mathcal{R}_{p}$. For patches of type $\mathcal{S}$,
$\mathcal{I}$ and $\mathcal{C}_2$ (Figure~\ref{fig:mappings}) the
canonical parameter space is taken to equal to the unit square
\begin{equation}\label{square}
  \mathcal{Q} \coloneqq [0,1]\times[0,1],
\end{equation}
whereas for patches of type $\mathcal{C}_1$ the L-shaped parameter
space
\begin{equation}\label{convex-param}
  \mathcal{L} \coloneqq ([0, 1] \times [\frac{1}{2}, 1]) \cup ([\frac{1}{2}, 1] \times [0, \frac{1}{2}])
\end{equation}
is used (Figure~\ref{fig:corner_mapping}). Simple Cartesian
discretizations of $\mathcal{Q}$ and $\mathcal{L}$ induce, via the
corresponding parametrizations used, the necessary discretizations of
each one of the patches $\Omega^{\mathcal{R}}_{p}$ for
$\mathcal{R}=\mathcal{S}$, $\mathcal{C}_1$, $\mathcal{C}_2$ and
$\mathcal{I}$, and for all relevant values of $p$.  The patch overlaps
are assumed to be ``sufficiently broad''---so that, roughly speaking,
except for discretization points near physical boundaries, each
discretization point $\mathbf{x}$, associated with a given patch
$\Omega^{\mathcal{R}}_p$, ($\mathcal{R}=\mathcal{S}$, $\mathcal{I}$,
$\mathcal{C}_1$, $\mathcal{C}_2$; $p=1,\dots,P_{\mathcal{R}}$) is
contained in a ``sufficiently interior'' region of at least one other
patch $\Omega^{\mathcal{R'}}_{p'}$ with
$(p',\mathcal{R}')\ne (p,\mathcal{R})$. This overlap condition is
imposed so as to ensure adequate interpolation ranges between patches
as well as creation of sufficiently smooth viscosity windows near
patch boundaries---as detailed in Sections~I--3.4
and~I--5.2. Section~\ref{subsec:overlap} below briefly reviews the
overlap requirement for all patches, and it presents additional
conditions required on the overlap between $\mathcal{S}$- and
$\mathcal{C}_1$-type patches.

In order to facilitate parallelization and to avoid the use of Fourier
expansions of extremely high order (so as to control differentiation
errors in the context of sharp flow features such as shock waves and
other solution discontinuities), the proposed algorithm allows for
partitioning of the patches $\Omega^{\mathcal{R}}_p$
($\mathcal{R}=\mathcal{S}$, $\mathcal{I}$, $\mathcal{C}_2$,
$\mathcal{C}_1$) into sets of similarly overlapping ``subpatches''
introduced in Section~\ref{subsec:overl_dec}.

The proposed domain-decomposition algorithm constructs at first
boundary-conforming patches
$\Omega^{\mathcal{R}}_{p}, p = 1, \dots, P_{\mathcal{R}}$, each one of
which is mapped via a domain mapping
\begin{equation}\label{eqn:mappings}
  \mathcal{M}^\mathcal{R}_{p}: \widetilde{\mathcal{R}} \rightarrow \overline{\Omega^{\mathcal{R}}_{p}}
\end{equation}
where, for $\mathcal{R}=\mathcal{S}$, $\mathcal{I}$ and
$\mathcal{C}_2$ (resp. $\mathcal{R}=\mathcal{C}_1$) we set
$\widetilde{\mathcal{R}} =\mathcal{Q}$ (resp.
$\widetilde{\mathcal{R}} =\mathcal{L}$); see
equations~\eqref{square} and~\eqref{convex-param}.

Using these notations and conventions the proposed general strategy
for construction of patch subdomains can be outlined as a sequence of
four steps.
\begin{enumerate}[(i)]
\item\label{pt1} $\mathcal{C}_1$-type patches: Patches around any existing
  $\mathcal{C}_1$ corner points are constructed following the
  procedure described in Section~\ref{subsec:C1-patches}. 
\item\label{pt2} $\mathcal{C}_2$-type patches: Patches around any
  existing $\mathcal{C}_2$ corner points are constructed following the
  procedure described in Section~\ref{subsec:C2-patches}.
\item\label{pt3} $\mathcal{S}$-type patches: Patches containing
  portions of smooth physical boundaries are then constructed,
  insuring a sufficiently wide overlap with the previously constructed
  $\mathcal{C}_1$-type and $\mathcal{C}_2$-type patches, as indicated
  in Section~I--3.1.
\item\label{pt4} $\mathcal{I}$-type patches: Interior patches, away from
  boundaries are finally constructed so as to cover the remaining
  interior of the domain $\Omega$, typically using affine mappings,
  and insuring a sufficiently wide overlap with existing
  $\mathcal{C}_1$-type, $\mathcal{C}_2$-type and $\mathcal{S}$-type patches
  introduced in steps~\eqref{pt1}, \eqref{pt2} and~\eqref{pt3}, as
  detailed in Section~I--3.1.
\item\label{pt5} Subpatches: The patches introduced in
  steps~\eqref{pt1}, \eqref{pt2}, \eqref{pt3} and~\eqref{pt4} are
  subsequently further decomposed into overlapping subpatches within
  each patch, as described in Section~\ref{L} for the patches
  described in step~\eqref{pt1}, and in Section~I--3.2.1 for
  the the patches described in steps~\eqref{pt2},~\eqref{pt3}
  and~\eqref{pt4}---if necessary, to ensure a sufficiently fine
  gridsize is achieved in each subpatch while using a constant number
  of discretization points per subpatch, as discussed in
  Sections~\ref{subsec:overlap} and \ref{subsubsec:subpatches}.
\end{enumerate}

\subsubsection{\label{subsec:C1-patches}$\mathcal{C}_1$- type patches}
\begin{figure}[H]
\centering
\includegraphics[width=0.8\linewidth]{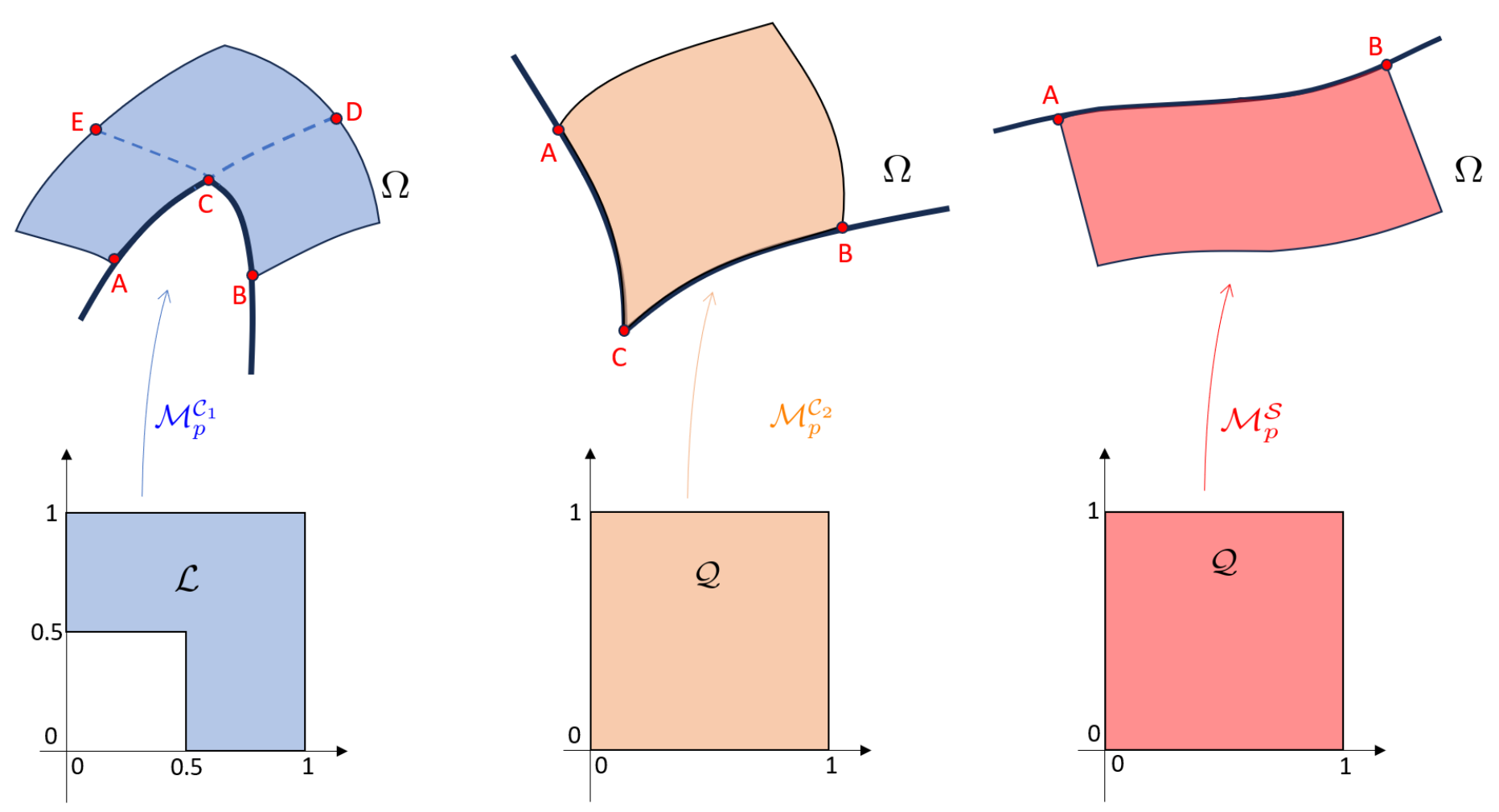}
\caption{Domain mappings for a $C_1$-type patch (left image pair), a
  $C_2$-type patch (middle image pair), and an $S$-type patch (right
  image pair).}
    \label{fig:corner_mapping}
\end{figure}

Let $C\in \Gamma$ denote a $\mathcal{C}_1$ corner point (see
Remark~\ref{CeqC2}). By definition, in a neighborhood of $C$ the
boundary $\Gamma$ can be represented as the union of two smooth arcs
$\overset{\frown}{AC}$ and $\overset{\frown}{BC}$ meeting at $C$, as
illustrated in the leftmost upper panel in
Figure~\ref{fig:corner_mapping}.  We assume, as we may, that smooth
and invertible 1D maps $\ell_A:[0,1]\to\overline{\Omega}$ and
$\ell_B:[0,1]\to\mathbb{R}^2$ are available that parametrize suitable
extensions of the curves $\overset{\frown}{AC}$ and
$\overset{\frown}{BC}$ beyond the point $C$, in such a way that, with
reference to Figure~\ref{fig:corner_mapping}, we have
\begin{equation}\label{corner_mapping-1}
\begin{split}
  \ell_A([0, \frac{1}{2}]) = \overset{\frown}{AC}, \qquad  \ell_A([\frac{1}{2}, 1]) = \overset{\frown}{CD}, \\
  \ell_B([0, \frac{1}{2}]) = \overset{\frown}{BC}, \qquad \ell_B([\frac{1}{2}, 1]) = \overset{\frown}{CE},
\end{split}
\end{equation}
(and, thus, in particular, $\ell_A([0, 1]) = \overset{\frown}{AD}$ and
$\ell_B([0, 1]) = \overset{\frown}{BE}$). Using these curves a 2D
parametrization
$\mathcal{M}^{\mathcal{C}_1}_p:\mathcal{L}\to\overline{\Omega^{\mathcal{C}_1}_p}$
of the form~\eqref{eqn:mappings} for a $\mathcal{C}_1$-type patch
$\Omega_p^{\mathcal{C}_1}$ around the point $C$, with parameter space as
depicted in the leftmost lower panel in
Figure~\ref{fig:corner_mapping}, and using the parameters
$\mathbf{q} = (q^1,q^2)$ is given by
\begin{equation}
  \label{eq:param_C2}
  \mathcal{M}^{\mathcal{C}_1}_p(q^1,q^2) = \ell_A(q^1)+\ell_B(q^2) - C.
\end{equation}

To conclude this section it is relevant to note that a straightforward
patching strategy around a $\mathcal{C}_1$ corner point obtained as a
union of two patches extruded separately from the curves
$\overset{\frown}{AD}$ and $\overset{\frown}{BE}$ results in a pair of
patches that cannot satisfy the inter-patch overlap condition required
per Section~\ref{subsec:overlap} below (that is  outlined in I--Remark~2 and illustrated in I--Figure~4; see also I--3.4.1)
---which, roughly speaking, requires that points in one patch that are
recipients of interpolation data from a second patch cannot themselves
be donors of interpolation data onto the second patch.

\subsubsection{\label{subsec:C2-patches}$\mathcal{C}_2$-type patches}

Let us now consider a $\mathcal{C}_2$ corner point $C\in \Gamma$ (see
Remark~\ref{CeqC2}). By definition, in a neighborhood of $C$ the
boundary $\Gamma$ can be represented as the union of two smooth arcs
$\overset{\frown}{AC}$ and $\overset{\frown}{BC}$ meeting at $C$,
which may be parametrized by smooth and invertible 1D maps
$\ell_A:[0,1]\to\mathbb{R}^2$ and $\ell_B:[0,1]\to\mathbb{R}^2$ in
such a way that, with reference to Figure~\ref{fig:corner_mapping}, we
have
\begin{equation}\label{corner_mapping}
  \ell_A([0, 1]) = \overset{\frown}{AC},\quad
  \ell_B([0, 1]) = \overset{\frown}{BC}.
\end{equation}
Using these parametrized curves we obtain the 2D parametrization
$\mathcal{M}^{\mathcal{C}_2}_p
:\mathcal{Q}\to\overline{\Omega^{\mathcal{C}_2}_p}$ given by
\begin{equation}
  \label{eq:param_C1}
  \mathcal{M}^{\mathcal{C}_2}_p(q^1,q^2) = \ell_A(q^1)+\ell_B(q^2) - C
\end{equation}
for a $\mathcal{C}_2$-type patch $\Omega_p^{\mathcal{C}_2}$ around the point
$C$.


\subsection{\label{subsec:overl_dec}Patches, subpatches and their grids}

The proposed algorithm is based on use of a decomposition of the
form~\eqref{decomp} into patches of various kinds, wherein, as in
Part~I, the patches used enjoy a ``sufficient amount of overlap''---in
the sense that every point $x\in\Omega^{\mathcal{R}}_p$
($\mathcal{R}=\mathcal{S}$, $\mathcal{C}_1$, $\mathcal{C}_2$, or
$\mathcal{I}$ and $1\leq p\leq P_{\mathcal{R}}$) that lies ``in the
vicinity'' of the boundary of $\Omega^{\mathcal{R}}_p$ must lie
``sufficiently deep'' within some patch $\Omega^{\mathcal{R}'}_{p'}$
in~\eqref{decomp} (i.e., within $\Omega^{\mathcal{R}'}_{p'}$ and away
from a vicinity of the boundary of $\Omega^{\mathcal{R}'}_{p'}$) for
some $(p',\mathcal{R}')\ne (p,\mathcal{R})$. The vicinity and depth
concepts alluded to above are defined and quantified in Part~I for
patches of type $\mathcal{C}$, $\mathcal{S}$ and $\mathcal{I}$ in
Sections~I--3.2.2, I--3.4.1 and~I--3.4.2---on the basis of the patch
discretizations used for such patches, which are themselves introduced
in Section~I--3.2.1. The discretization and depth concepts used in
this paper for patches of types $\mathcal{C}_2$, $\mathcal{S}$ and
$\mathcal{I}$ coincide with the ones used in Part~I for patches of
types $\mathcal{C}$, $\mathcal{S}$ and $\mathcal{I}$,
respectively. The corresponding concepts for patches of type
$\mathcal{C}_1$, in turn, are introduced in
Section~\ref{subsec:overlap} on the basis of the corresponding
discretization structures introduced in Section~\ref{L}. Roughly
speaking, both the discretizations used in Part~I and in the present
paper are produced by first introducing, at the level of the parameter
space, a set of overlapping subpatches of each patch, and the depth
concept is defined in terms of certain layers of discretization points
near patch boundaries. (The condition of overlap between patches of
type $\mathcal{C}_1$ and other patches includes an additional
requirement to ensure stability of the overall time-stepping
algorithm---namely, that, as discussed in
Section~\ref{subsec:overlap}, the overlapping $\mathcal{S}$-type
patches must extend up to and including the corner point.)

As in Part~I, for flexibility in both geometrical representation and
discretization refinement, the patching structure we use incorporates
a decomposition in ``sub-patches'' $\Omega^{\mathcal{R}}_{p, \ell}$ of
each patch $\Omega^{\mathcal{R}}_{p}$ ($\mathcal{R}=\mathcal{S}$,
$\mathcal{C}_1$, $\mathcal{C}_2$ or $\mathcal{I}$):
\[
  \Omega^{\mathcal{R}}_{p} = \bigcup_{\ell=1}^{r_p} \Omega^{\mathcal{R}}_{p, \ell},
\]
on each of which a patch-dependent number of discretization points is
used. The sub-patching strategy used in the present paper for patches
of types, $\mathcal{C}_2$, $\mathcal{S}$ and $\mathcal{I}$ coincides
with the ones used in Part~I for the $\mathcal{C}$, $\mathcal{S}$ and
$\mathcal{I}$ patch types, and the corresponding strategy for
$\mathcal{C}_1$-type patches, which is described in Section~\ref{L}
below, is a natural variant thereof. Briefly, in all cases subpatching
is set up at the parameter-space level, and it is achieved by
partitioning the canonical parameter spaces $\mathcal{P}= \mathcal{Q}$
and $\mathcal{P}= \mathcal{L}$, as described in Sections~\ref{Q}
and~\ref{L} below, so as to ensure sufficiently large amounts of overlap
between neighboring subpatches.

\subsubsection{\label{Q}$\mathcal{Q}$ parameter space discretization
(for $\mathcal{S}$-, $\mathcal{C}_2$- and $\mathcal{I}$-type patches) }

Discretizations for each one of the patches $\Omega^\mathcal{R}_p$
($\mathcal{R}=\mathcal{S}$, $\mathcal{C}_2$, $\mathcal{I}$;
$1\leq p\leq P_\mathcal{R}$) are obtained by introducing suitable
grids in the parameter space $\mathcal{Q}$ (equation~\eqref{square}),
as indicated in Section~I--3.2.1.

\subsubsection{\label{L}$\mathcal{L}$ parameter space discretization
(for $\mathcal{C}_1$-type patches) }
In analogy to the approach utilized in Section~\ref{Q}, in order to
obtain discretizations for each one of the patches
$\Omega^{\mathcal{C}_1}_p$ ; $1\leq p\leq P_{\mathcal{C}_1}$, we
introduce the parameter space grids
\begin{equation}
  \label{eq:param-grid-C1}
  G^{\mathcal{C}_1}_{p} = \mathcal{L} \cap \big\{(q^{\mathcal{C}_1, 1}_{p, i}, q^{\mathcal{C}_1, 2}_{p, j})\ :\ q^{\mathcal{C}_1, 1}_{p, i} = ih^{\mathcal{C}_1, 1}_{p}, q^{\mathcal{C}_1, 2}_{p, j} =
  jh^{\mathcal{C}_1, 2}_p, \quad (i, j) \in \mathbb{N}_0^2 \big\}
\end{equation}
(where $ \mathbb{N}_0 = \mathbb{N}\cup\{ 0\}$) which, in the present
case, we take to be {\em equiaxed}---that is, with equal parameter
spacegrid sizes along the $q^1$ and $q^2$ directions, as detailed in
Remark~\ref{square_subpatches}.
\begin{remark}\label{square_subpatches}
  To enable use of a single $\mathcal{L}$ parameter space for all
  $\mathcal{C}_1$-type patches, as in the approach used for the
  $\mathcal{Q}$ parameter space described in Part~I, it is assumed
  here that each $\mathcal{C}_1$-type patch in physical space is
  constructed in such a way that its two arms are approximately equal
  in width and length. Under this assumption it is reasonable to
  utilize a single grid size parameter
  \begin{equation}\label{eq:eq-mesh-L}
  h^{\mathcal{C}_1}_{p} = h^{\mathcal{C}_1, 1}_{p} = h^{\mathcal{C}_1, 2}_{p} = \frac{1}{N^{\mathcal{C}_1}_p + 2n_v - 1},
\end{equation}
for $\mathcal{C}_1$-type patches, where the integers
$N^{\mathcal{C}_1}_p$ and $n_v$ are described in what follows. The
pair of meshsize parameters $h^{\mathcal{C}_1, 1}_{p}$ and
$h^{\mathcal{C}_1, 2}_{p}$, which is analogous to the pairs of
parameters $h^{\mathcal{R}, 1}_{p}$ and $h^{\mathcal{R}, 2}_{p}$ with
$\mathcal{R}=\mathcal{S}$, $\mathcal{C}$ or $\mathcal{I}$ introduced
in Section~I--3.2.1, is preserved here as well, however---in spite of
the fact that, per equation~\eqref{eq:eq-mesh-L}, these quantities are
taken to coincide---in order to streamline the notation in connection
with the viscosity assignment strategy outlined in
Section~\ref{sec:viscosity} and detailed in Section~I--5.
\end{remark}

\begin{figure}[H]
  \centering
  \includegraphics[width=0.7\linewidth,]{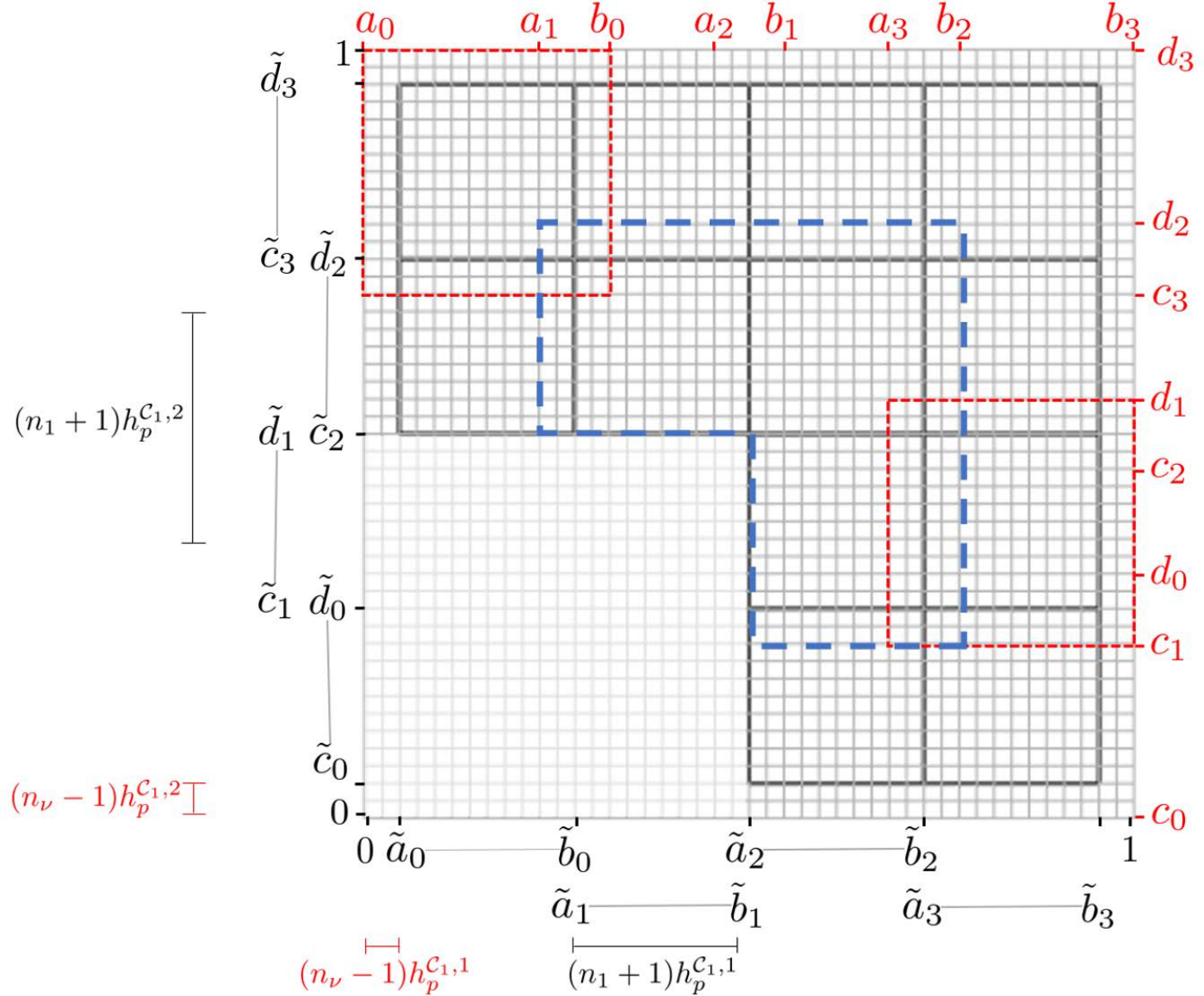}
  \caption{Subpatch decomposition of the $\mathcal{L}$-parameter-space
    used for patches $\Omega^{\mathcal{R}}_p$ with
    $\mathcal{R} = \mathcal{C}_1$; for this illustration the geometric
    parameter values $r^{\mathcal{C}_1}_p = 4$, $n_1 = 9$, and
    $n_v = 3$ were used. According to the text prescriptions the patch
    is partitioned into a total of
    $\frac{3}{4} \big( r^{\mathcal{C}_1}_p \big)^2 = 12$
    non-overlapping preliminary square subpatches (with boundaries
    shown as black solid lines). A total of 10 overlapping patches
    exist for this configuration, three of which are displayed,
    namely, two overlapping patches (shown in red) separated from the
    corner point and the physical boundary
    (cf. Figure~\ref{fig:corner_mapping} left), and the corner patch
    $H_0 = L^{\mathcal{C}_1}_p$, depicted as a blue dashed
    polygon. For the preliminary subpatches adjacent to the physical
    boundary but not adjacent to the corner (2 of which exist in this
    example), the associated overlapping subpatches (not displayed)
    are rectangular in shape as they are not allowed to penetrate the
    obstacle.  Discretization lines are shown in gray.}
  \label{fig:lpatch}
\end{figure}

The boundary vicinity concept and associated parameter $n_v$
considered here are closely related to the analogous concepts for
patches of types $\mathcal{R}=\mathcal{S}$, $\mathcal{C}$,
$\mathcal{I}$ introduced in Section~I--3.2.1, whose parameter space is
the square $\mathcal{Q}$ and wherein the boundary vicinity equals the
union of four thin rectangles. As in Part~I, the value $n_v=9$ is used
throughout this paper for all patch types considered, i.e., for
patches of type $\mathcal{R}=\mathcal{S}$, $\mathcal{C}_1$
$\mathcal{C}_2$ and $\mathcal{I}$. With reference to
equation~\eqref{convex-param} and Figure~\ref{fig:lpatch}, on the
other hand, the parameter $N^{\mathcal{C}_1}_p$ in
equation~\eqref{eq:eq-mesh-L} is defined so that the complete interval
$[0,1]$ in either the $q^1$ or $q^2$ directions contains
$N^{\mathcal{C}_1}_p + 2 n_v$ discretization points.

The strategy used for sub-patch decomposition of $\mathcal{C}_1$-type
patches, in turn, is based on a decomposition of the corresponding
parameter space $\mathcal{P} = \mathcal{L}$, as illustrated in
Figure~\ref{fig:lpatch} and described in what follows. Following the
lines of the $\mathcal{P} = \mathcal{Q}$ decomposition presented in
Section~I--3.2.1, here we utilize, for a given even integer
$r^{\mathcal{C}_1}_p$, a union of non-overlapping squares, called {\em
  preliminary subpatches}, strictly contained within the canonical
L-shaped domain $\mathcal{L}$:

\begin{equation}
  \label{eq:squares}
 \bigcup_{(r, s) \in \widetilde{\Theta}^{\mathcal{C}_1}_{p}} [\widetilde{a}_{r}, \widetilde{b}_{r}] \times [\widetilde{c}_{s}, \widetilde{d}_{s}]   \subsetneq \mathcal{L},\qquad\qquad (\widetilde{\Theta}^{\mathcal{C}_1}_{p} = \widetilde{\Theta}^{\mathcal{C}_1,
  1}_{p}\cup \widetilde{\Theta}^{\mathcal{C}_1, 2}_p\cup
\widetilde{\Theta}^{\mathcal{C}_1, 3}_{p}).
\end{equation}
Here, using the integer parameters $r^{\mathcal{C}_1}_p$ and
$s^{\mathcal{C}_1}_p$ which in the present case
$\mathcal{R}=\mathcal{C}_1$ (and in contrast with $r^{\mathcal{R}}_p$
and $s^{\mathcal{R}}_p$ parameters used in the cases
$\mathcal{R}=\mathcal{S}$, $\mathcal{C}_2$ and $\mathcal{I}$) are
taken to be {\em equal},
\begin{equation}
  \label{eq:equality}
  r^{\mathcal{C}_1}_p = s^{\mathcal{C}_1}_p
\end{equation}
(see Remark~\ref{square_subpatches}), we have set
\begin{equation}\label{eq:sp_index2}
  \begin{aligned}
    &\widetilde{\Theta}^{\mathcal{C}_1, 1}_{p} = \{(r,s)\in\mathbb{N}_0^2\ |\ 0\leq r\leq \frac{r^{\mathcal{C}_1}_p}{2}
      - 1 \mbox{ and } \frac{s^{\mathcal{C}_1}_p}{2} \leq  s \leq s^{\mathcal{C}_1}_p - 1 \}, \\
    &\widetilde{\Theta}^{\mathcal{C}_1, 2}_{p} = \{(r,s)\in\mathbb{N}_0^2\ |\ \frac{r^{\mathcal{C}_1}_p}{2} \leq r\leq r^{\mathcal{C}_1}_p
      - 1 \mbox{ and } \frac{s^{\mathcal{C}_1}_p}{2} \leq s \leq s^{\mathcal{C}_1}_p - 1 \}, \\
    &\widetilde{\Theta}^{\mathcal{C}_1, 3}_{p} = \{(r,s)\in\mathbb{N}_0^2\ |\ \frac{r^{\mathcal{C}_1}_p}{2} \leq r\leq r^{\mathcal{C}_1}_p
      - 1 \mbox{ and } 0\leq s\leq \frac{s^{\mathcal{C}_1}_p}{2} - 1 \};
  \end{aligned}
\end{equation}
clearly, $r^{\mathcal{C}_1}_p = s^{\mathcal{C}_1}_p$ is the number of
preliminary (non-overlapping) subpatches along each of the arms of the
L-shaped parameter space. The number of discretization points along
each side of the (square) preliminary subpatches will be denoted by
$n_1$---so that, with reference to~\eqref{eq:squares}, for all $r$ and
$s$ we have
$(\widetilde{b}_{r}-\widetilde{a}_{r}) =
(\widetilde{d}_{s}-\widetilde{c}_{s}) =(n_1 +
1)h^{\mathcal{C}_1}_{p}$. Finally, letting
$N^{\mathcal{C}_1}_p = r^{\mathcal{C}_1}_p(n_1 + 1) - 1$ the total
number of discretization points along each arm of the L-shaped
parameter space equals $N^{\mathcal{C}_1}_p+2n_v$---as indicated above in this section.

In view of the setup just described it is easy to check that the
following relations hold:
\begin{equation}\label{discretization_L}
  \begin{aligned}
   & \tilde{a}_0 = (n_{v}-1) h^{\mathcal{C}_1, 1}_p, &\\
   & \tilde{b}_{r} - \tilde{a}_r = (n_1 + 1) h^{\mathcal{C}_1, 1}_p, &\qquad (0 \leq r \leq r^{\mathcal{C}_1}_p - 1),  \\
   & \tilde{a}_{r+1} = \tilde{b}_{r},  &\qquad (0 \leq r < r^{\mathcal{C}_1}_p - 1),  \\
   & \tilde{c}_0 = (n_{v}-1) h^{\mathcal{C}_1, 2}_p, & \\
   & \tilde{d}_{s} - \tilde{c}_s = (n_1+1)  h^{\mathcal{C}_1, 2}_p,  &\qquad (0 \leq s \leq s^{\mathcal{C}_1}_p - 1 ), \\
   & \tilde{c}_{s+1} = \tilde{d}_{s},  &\qquad (0 \leq s < s^{\mathcal{C}_1}_p - 1).
  \end{aligned}
\end{equation}

Typically the selection $n_1\approx \frac 12 n_0$ is utilized, where,
per Section~I--3.2.1, $n_0$ equals the number of discretization points
along each side of the (possibly non-square, non-overlapping)
preliminary subpatches for patches of type $\mathcal{R}=\mathcal{S}$,
$\mathcal{I}$ and $\mathcal{C}_2$ (noting once again that the
$\mathcal{C}_2$ patch type is called the $\mathcal{C}$ patch type in
Part~I); the values $n_0 = 83$ (also used in Part~I) and
$n_1 = 43\approx n_0/2$ are used in all the test cases considered in
this paper.

As for the previously considered patch types, the {\em overlapping}
subpatches associated with a given $\mathcal{C}_1$-type patch are
characterized, via the patch parametrization, by an overlapping
decomposition of $\mathcal{L}$ as the union of closed subpatches. In
the present case, in addition to overlapping rectangular subpatches,
as are used in connection with patches of type
$\mathcal{R}=\mathcal{S}$, $\mathcal{I}$ and $\mathcal{C}_2$, an
``L-shaped'' subpatch $L^{\mathcal{C}_1}_p$ is utilized. In detail,
the overlapping sub-patch decomposition of a $\mathcal{C}_1$ patch is
constructed on the basis of the parameter-space decomposition
\begin{equation}
  \label{eq:rect_L}
  \mathcal{L} =  
    \bigcup_{q \in \Theta^{\mathcal{C}_1}_p} H_q, \qquad (\Theta^{\mathcal{C}_1}_p = \Theta^{\mathcal{C}_1, 1}_{p} \cup \Theta^{\mathcal{C}_1, 2}_{p} \cup \Theta^{\mathcal{C}_1, 3}_{p} \cup \{0\})
\end{equation}
where
\begin{equation}\label{H}
  \begin{alignedat}{3}
    &H_q = [a_r, b_r] \times [c_s, d_s]  &\quad  &\mbox{for} &\quad  &  q = (r, s) \in  \Theta^{\mathcal{C}_1, 1}_{p} \cup \Theta^{\mathcal{C}_1, 2}_{p} \cup \Theta^{\mathcal{C}_1, 3}_{p} \\
    &H_q = L^{\mathcal{C}_1}_p &  &\mbox{for} &  & q = 0,  
   \end{alignedat}
 \end{equation}
where the ($p$-dependent) endpoints $a_{r}$, $b_{r}$, $c_{s}$ and
$d_{s}$ are given by
\begin{equation}\label{subpatchesQ}
  \begin{aligned}
    &a_{r}  = \tilde{a}_r - (n_{v}-1) h^{\mathcal{C}_1}_p, &\qquad (0 \leq r \leq r^{\mathcal{C}_1}_p - 1),  \\
    &b_{r}  = \tilde{b}_{r} + (n_{v}-1) h^{\mathcal{C}_1}_p, &\qquad (0 \leq r \leq r^{\mathcal{C}_1}_p - 1), \\
    &c_{s}  = \tilde{c}_s - (n_{v}-1) h^{\mathcal{C}_1}_p, &\qquad (0 \leq s \leq s^{\mathcal{C}_1}_p - 1),  \\
    &d_{s}  = \tilde{d}_{s} + (n_{v}-1) h^{\mathcal{C}_1}_p, &\qquad (0 \leq s \leq s^{\mathcal{C}_1}_p - 1),  
  \end{aligned}
\end{equation}
and where the index sets
\begin{equation}\label{eq:ind_sub}
  \begin{aligned}
  & \Theta^{\mathcal{C}_1, 1}_{p} = \widetilde{\Theta}^{\mathcal{C}_1, 1}_{p} \setminus  \big\{(\frac{r^{\mathcal{C}_1}_p}{2}
  - 1, \frac{s^{\mathcal{C}_1}_p}{2})\big\}, \\
  & \Theta^{\mathcal{C}_1, 2}_{p} = \widetilde{\Theta}^{\mathcal{C}_1, 2}_{p} \setminus  \big\{(\frac{r^{\mathcal{C}_1}_p}{2} , \frac{s^{\mathcal{C}_1}_p}{2})\big\}, \\
  & \Theta^{\mathcal{C}_1, 3}_{p} = \widetilde{\Theta}^{\mathcal{C}_1, 3}_{p} \setminus  \big\{(\frac{r^{\mathcal{C}_1}_p}{2} , \frac{s^{\mathcal{C}_1}_p}{2} - 1)\big\},
  \end{aligned}
\end{equation}
were used. The corner set
$H_0 = L^{\mathcal{C}_1}_p$ in~\eqref{H}, finally, is defined by

\begin{equation} \label{eq:inner_L}
    L^{\mathcal{C}_1}_p = L^{\mathcal{C}_1}_{p, 1} \cup L^{\mathcal{C}_1}_{p, 2} \cup L^{\mathcal{C}_1}_{p, 3},
\end{equation}
where, calling temporarily $m_p = \frac{r^{\mathcal{C}_1}_p}{2}$ to
improve readability, we have called
\begin{equation} \label{eq:inner_L_dec}
  \begin{aligned}
   & L^{\mathcal{C}_1}_{p, 1} = [a_{m_p - 1},
  b_{m_p - 1}] \times
  [c_{m_p},
  d_{m_p}], \\
   & L^{\mathcal{C}_1}_{p, 2} = [a_{m_p},
  b_{m_p}] \times
  [c_{m_p - 1}, d_{m_p
    - 1}], \\
   & L^{\mathcal{C}_1}_{p, 2} = [a_{m_p - 1},
  b_{m_p - 1}] \times
  [c_{m_p},
  d_{m_p}].
  \end{aligned}
\end{equation}
Note the three singleton sets subtracted in~\eqref{eq:ind_sub} contain
the three index pairs $q$ corresponding to the squares $H_q$ whose
union contains the corner subpatch $H_0 = L^{\mathcal{C}_1}_p$.

For $\ell \in \Theta^{\mathcal{C}_1}_p$, we denote by
$G^{\mathcal{C}_1}_{p, \ell}$ the parameter space grid for the $\ell$-th
subpatch of the patch $\Omega^{\mathcal{C}_1}_p$:
\[
  G^{\mathcal{C}_1}_{p, \ell} =   G^{\mathcal{C}_1}_{p} \cap H_{\ell}.
\]
Using the mapping $\mathcal{M}^{\mathcal{C}_1}_p$ described in
Section~\ref{subsec:C1-patches} we obtain the desired subpatches
$\Omega^{\mathcal{C}_1}_{p, \ell} = \mathcal{M}^{\mathcal{C}_1}_p(
H_{\ell})$ as well as the physical grids
\begin{equation}
  \label{eq:c1_grid}
    \mathcal{G}^{\mathcal{C}_1}_p = \mathcal{M}^{\mathcal{C}_1}_p\left( G^{\mathcal{C}_1}_p \right) \quad \mbox{and} \quad \mathcal{G}^{\mathcal{C}_1}_{p, \ell} = \mathcal{M}^{\mathcal{C}_1}_p\left( G^{\mathcal{C}_1}_{p, \ell}\right);
\end{equation}
note that the set
\begin{equation} \label{D_C1}
\mathcal{D}^{\mathcal{C}_1}_{p, \ell} = \{(i, j) \in \mathbb{N}_0^2 \ |\ (ih^{\mathcal{C}_1, 1}_p, jh^{\mathcal{C}_1, 2}_p) \in G^{\mathcal{C}_1}_{p, \ell} \}
\end{equation}
equals the set of indices of all grid points in the set
$\mathcal{G}^{\mathcal{C}_1}_{p, \ell}$.

\subsubsection{\label{subsec:overlap}Minimum subpatch overlap
  condition and special $\mathcal{S}$-$\mathcal{C}_1$ patch overlaps}

As detailed in Sections~I--3.4 and~I--5.2, and as alluded to
previously in Section~\ref{subsec:decomposition}, in order to properly
enable inter-patch data communication and multi-patch viscosity
assignment propagation, the overlap between pairs of neighboring
patches $\Omega^{\mathcal{R}}_p$ and $\Omega^{\mathcal{R'}}_{p'}$
($(p',\mathcal{R}')\neq (p,\mathcal{R})$) must be sufficiently
broad. In the cases $\mathcal{R},\mathcal{R}'\neq\mathcal{C}_1$ the
description coincides with the one presented in Section~I--3.2.2
($\mathcal{C}_2$ is called $\mathcal{C}$ in Part~I), while the case
$\mathcal{R}=\mathcal{C}_1$ requires additional considerations. As in
Section~I--3.2.2, in all cases, including the case
$\mathcal{R}=\mathcal{C}_1$, the required overlap breadth is
quantified in terms of a parameter associated with the selected
subpatches of each patch (namely, an integer parameter $n_v$
introduced in Section~I--3.2.1 and in Section~\ref{L} above. The use
of the parameter $n_v$ and the concept of patch-boundary ``sides''
introduced in Section~I--3.2.2 to characterize minimum subpatch
overlaps extend directly to the context of
$\mathcal{R}= \mathcal{C}_1$-type patches introduced in
Section~\ref{L}. The ``sides'' of a patch (resp. of a subpatch) are
defined as the images under the patch parametrization of each one of
the straight segments that make up the boundary of the corresponding
parameter polygon $\mathcal{Q}$ or $\mathcal{L}$,
cf. Remark~\ref{rem_bdries} (resp. each one of the sides of the
rectangles in equations~I--(21) and~\eqref{H} and the sides of the
L-shaped polygons~\eqref{eq:inner_L}). Utilizing these concepts, the
minimum-overlap condition is said to be satisfied for a given patch
$\Omega^{\mathcal{R}}_p$ if and only if for each side of
$\Omega^{\mathcal{R}}_p$ that is not contained in $\partial \Omega$
(cf. Remark~\ref{rem_bdries}), a $(2 n_v + 1)$-point wide layer of
grid points adjacent to that side in the $\Omega^{\mathcal{R}}_p$ grid
is also included in a union of one or more patches
$\Omega^{\mathcal{R'}}_{p'}$ with
$(p',\mathcal{R}')\neq (p,\mathcal{R})$. (Note that, per the
construction in Sections~\ref{Q} and~\ref{L} at the parameter space
level, subpatches of a single patch satisfy the minimum overlap
condition embodied in equations~\eqref{subpatchesQ}
and~\eqref{discretization_L}: neighboring subpatches share a
$(2 n_v + 1)$-point wide overlap.)

A particular note is necessary in regard to the overlap of
$\mathcal{C}_1$- and $\mathcal{S}$-type patches. Numerical tests have
shown that spurious oscillations and instability may emanate from
boundary regions at the intersection of a $\mathcal{C}_1$ and an
$\mathcal{S}$ patch near a corner point. Fortunately, however,
numerical experimentation for the wedge and prism geometries
(Figures~\ref{Wedge_Flow_solutions} and~\ref{Prism_Flow_solutions})
suggests that such spurious oscillations may be eliminated provided
the $\mathcal{S}$-type patch extends along the boundary up to and
including the corner point (e.g. including the complete arcs
$\overset{\frown}{AC}$ and $\overset{\frown}{BC}$ on the upper-left
panel of Figure~\ref{fig:corner_mapping}).  With such a prescription
the observed spurious oscillations are eliminated and the stability of
the algorithm is recovered.

\subsubsection{Subpatches and grid refinement\label{subsubsec:subpatches}}

In analogy with the discussion presented in Section~I--3.2.3, here we
note that general overlapping patch decompositions satisfying the
minimum overlap condition introduced in Section~\ref{subsec:overlap}
(which itself extends Section~I--3.2.2 to the context of domains with
corners) can be obtained on the basis of the procedures described in
Section~\ref{subsec:decomposition}. As briefly discussed below in this
section (which closely parallels Section~I--3.2.3 while adding details
applicable to $\mathcal{C}_1$-type patches), further, the strategies
presented in Sections~\ref{Q} and~\ref{L} can be used to produce sets
of overlapping subpatches, each one endowed with a subpatch grid, in
such a way that, 1)~The set of all patches (namely, all
$\mathcal{S}$-, $\mathcal{C}_1$-, $\mathcal{C}_2$- and
$\mathcal{I}$-type patches) satisfies the overlap conditions
introduced in Section~\ref{subsec:overlap}; 2)~Each subpatch of a
$\mathcal{C}_1$-type patch (resp. of an $\mathcal{R}$-type patch with
$\mathcal{R}=\mathcal{S}$, $\mathcal{C}_2$ and $\mathcal{I}$),
contains no more than (resp. exactly)
\begin{equation}
  \label{eq:NQ}
  N_{\mathcal{Q}}=\left(n_0 + 2n_v\right)^2
\end{equation}
discretization points; and, 3)~All of the subpatch parameter space
grid sizes $h^{\mathcal{R}, 1}_p$ and $h^{\mathcal{R}, 2}_p$
($\mathcal{R}=\mathcal{S}$, $\mathcal{C}_1$, $\mathcal{C}_2$,
$\mathcal{I}$; $1\leq p\leq P_\mathcal{R}$) are such that the
resulting physical grid sizes (defined as the maximum distance between
two consecutive grid points in physical space) are less than or equal
to a user-prescribed upper bound $\overline{h}>0$.

As indicated in Section~\ref{subsec:overlap}, the sets of patches and
subpatches constructed per the descriptions in Sections~\ref{Q}
and~\ref{L} satisfy points~1) and~2) above, but, clearly, they may or
may not satisfy point~3) for certain values of $\mathcal{R}$ and
$p$. For such $(\mathcal{R},p)$ patches the integers
$r^{\mathcal{R}}_p$ and $ s^{\mathcal{R}}_p$ (cf.  Sections I-3.2.1
and~\ref{L}) are suitably increased---thus proportionally increasing
the numbers of subpatches and the number of gridpoints of the patch
$\Omega^{\mathcal{R}}_p$ while leaving the overall patch decomposition
and mappings $\mathcal{M}^{\mathcal{R}}_p$ unchanged---until the
physical grid size upper bound condition in point~3) is satisfied, as
desired. Additional grid refinements can be obtained if needed, as
indicated in the last paragraph of Section~I--3.2.3.

\subsection{\label{subsec:connectivity} Patch/subpatch communication
  of solution values}

The overlap built into the patch/subpatch domain decomposition
approach described in the previous sections enables the communication
(possibly via interpolation) of subpatch grid-point values of the
solution $\textbf{e}$ in the vicinity of boundaries of
subpatches---which is a necessary element in the Runge-Kutta based
overlapping-patch time-stepping algorithm we use (cf. Section~I--4).
The communication of solution values relies on the concept of
``$n_f$-point fringe region'' of a subpatch, which was defined in
Section~I--3.4 for patches of types $\mathcal{S}$, $\mathcal{C}_2$ and
$\mathcal{I}$, and which is extended here to include patches of type
$\mathcal{C}_1$.  To introduce this concept we consider the sides of
any subpatch $\Omega^{\mathcal{R}}_{p, \ell}$ (defined in
Section~\ref{subsec:overlap}), and we note that, per the constructions
in that section, each side of a subpatch is either contained within
$\Gamma$ or it intersects $\Gamma$ at most at one point. In the first
(resp. the second) case we say that the side is ``external''
(resp. ``internal''). Notably for a $\mathcal{C}_1$-type patch, which,
per the steps outlined at the end of
Section~\ref{subsec:decomposition}, discretizes the vicinity of a
$\mathcal{C}_1$ corner point of ${\Omega}$, (resp. a subpatch of a
$\mathcal{C}_1$-type patch), the sides equal to the arcs
$\overset{\frown}{AC}$ and $\overset{\frown}{BC}$ (resp. the sides
included in the arcs $\overset{\frown}{AC}$ and
$\overset{\frown}{BC}$) (see upper-left panel of
Figure~\ref{fig:corner_mapping}) are always considered to be external.

For each subpatch $\Omega^{\mathcal{R}}_{p, \ell}$
($\mathcal{R}=\mathcal{S}$, $\mathcal{C}_1$, $\mathcal{C}_2$,
$\mathcal{I}$; $1\leq p\leq P_\mathcal{R}$;
$\ell \in \Theta^{\mathcal{R}}_p$), we define the $n_f$-point fringe
region $ \mathcal{F}^{\mathcal{R}}_{p, \ell, n_f}$ analogously to
Section~I--3.4.  Roughly speaking, the $n_f$-point fringe region of
$\Omega^{\mathcal{R}}_{p, \ell}$ is a set of discretization points
$(x_i,y_j)\in \Omega^{\mathcal{R}}_{p, \ell}$ that are no more than
$n_f$-points away from an internal boundary of the subpatch---in the
sense that, starting from $(x_i,y_j)$ there exists a sequence of
$n_f-1$ additional consecutive discretization points along one of the
two parameter space dimensions, the last one of which is contained in
an internal subpatch boundary. As indicated in Sections I--3.4,
I--3.4.1, I--3.4.2 and Remark~2 in Part~I, it is important to take
into account the following considerations~1)~All the implementations
presented in the two-paper sequence use the value $n_f = 5$; 2)~The
time-stepping method utilizes two different algorithms for the
communication of solution data between pairs of subpatches, including
``inter-patch'' communication (that is, communication between two
subpatches of different patches, which requires polynomial {\em
  interpolation} of grid-point values of the solution $\textbf{e}$),
and ``intra-patch'' communication (which only involves {\em exchange}
of grid-point values of $\textbf{e}$ between two subpatches of the
same patch); and, 3)~The satisfaction of the minimum overlap condition
(Section~\ref{subsec:overlap}) ensures that data donor grid points are
not themselves receivers of data from other patches. (The stability of
the time-stepping method would be compromised if donor points were
also recipients of data from other patches.)
\begin{algorithm}
  \begin{algorithmic}[1]
    \Begin
    \State \textbackslash \textbackslash Initialization.
    \State Execute Algorithm~\ref{alg_init} to initialize the mesh and the multi-patch FC-SDNN solver. 
    \State Initialize time: $t = 0$.
     \State \textbackslash \textbackslash Time stepping.
  \While {$t < T$}  
  \State\label{pre_vis} Execute Algorithm~\ref{alg_prelim_artvisc} (parallel evaluation of the preliminary artificial viscosity operator  $\widehat\mu[\textbf{e}^{\mathcal{R}}_{p, \ell}]$).
  \State\label{visc_alg} Communicate the preliminary artificial viscosity values (line~\ref{pre_vis}) to all patches and subpatches and evaluate the overall artificial viscosity  $\mu^{\mathcal{R}}_{p,
    \ell}[\textbf{e}]$, in parallel, in each subpatch grid (see point~\ref{grids} in Section~\ref{main-elms}) as required by Algorithm~\ref{alg_stage}
  \State Execute Algorithm~\ref{alg_filt} (parallel filtering of the solution vector $\textbf{e}^{\mathcal{R}}_{h, p, \ell}$).
  \State Evaluate the temporal step-size $\Delta t$ (see point~\ref{time-stepping} in Section~\ref{main-elms}).
  \For {each stage of the SSPRK-4 time step}
  \State Execute Algorithm~\ref{alg_stage} (parallel evaluation of the the  SSPRK-4 stage  on the solution vector  $\textbf{e}^{\mathcal{R}}_{h, p, \ell}$  and enforcement of boundary conditions).
  \State\label{sol_comm} Communicate the solution values between neighboring patches and subpatches via exchange and interpolation, as relevant (see Section~\ref{subsec:connectivity}).
  \EndFor
  \State Update time: $t = t + \Delta t$
  \State Write solution values to disk at specified time steps $t$.
  \EndWhile
  \End
  \caption{Overall FC-SDNN driver program}
\label{alg_total}
\end{algorithmic}
\end{algorithm}

\section{\label{sec:viscosity}Multi-patch artificial viscosity assignment}

The multi-patch artificial viscosity-assignment algorithm used in this
paper is essentially identical to that presented in Section~I--5: the
description presented in that section is applicable in the present
context provided the the patch-type list ``$\mathcal{S}$,
$\mathcal{C}$, $\mathcal{I}$'' in Section~I--5 is replaced by the
corresponding patch list ``$\mathcal{S}$, $\mathcal{C}_1$,
$\mathcal{C}_2$, $\mathcal{I}$'' used here. Roughly speaking, the
viscosity-assignment method used in both Part~I and the present
Part~II papers obtains the artificial viscosity throughout $\Omega$ by
first computing, for each subpatch $\Omega^{\mathcal{R}}_{p, \ell} $,
a preliminary subpatch-wise viscosity assignment at each
discretization point located in the interior of the subpatch (or more
precisely, in the set
$\mathcal{G}^{\mathcal{R}}_{p, \ell} \setminus
\mathcal{F}^{\mathcal{R}}_{p, \ell, n_f}$ consisting of points located
outside of the $n_f$-point fringe region of the subpatch, as described
in Section~\ref{subsec:connectivity}). This preliminary
single-subpatch viscosity-assignment method, which relies on an
Artificial Neural Network (ANN) for classification of the smoothness
of the solution, is essentially identical to the single-patch
procedure first introduced in~\cite{bruno2022fc} (which, in turn, is
based on a number of the important ANN elements introduced
in~\cite{schwander2021controlling} in the context of periodic
solutions). The preliminary viscosity assignments obtained on each
patch are then combined across all patches by means of a
smoothing-blending operator $\Lambda$ constructed on the basis of a
set of multi-patch windowing functions---a procedure which uniquely
defines the artificial viscosity throughout $\Omega$.

\section{Overview of the multi-patch  FC-SDNN algorithm\label{main-elms}}

The FC-SDNN method relies on a number of elements that are briefly
reviewed below---to provide a summary of the method as well as a
quick-reference source for it.
\begin{enumerate}
\item\label{grids} The overall FC-SDNN algorithm relies on a
  decomposition of the computational domain into a number of
  overlapping patches $\Omega^{\mathcal{R}}_{p}$ with
  $\mathcal{R} = \mathcal{S}$, $\mathcal{C}_1$, $\mathcal{C}_2$,
  $\mathcal{I}$; $p=1,\dots,P_{\mathcal{R}}$, (including interior and
  smooth boundary patches considered in Section~I--3.2 as well as
  corner patches considered in Section~\ref{subsec:overl_dec}), as
  well as associated patch parametrizations
  $\mathcal{M}^{\mathcal{R}}_{p}$, and Jacobians
  $J_{\mathbf{q}}^{\mathcal{R}, p}$ of the inverse mappings
  $\big(\mathcal{M}^{\mathcal{R}}_{p}\big)^{-1}$, that are used to
  express the Euler equations in parameter space, as detailed in
  (Section~I--3.3). Overlapping subpatches of these patches are then
  constructed, and discretizations
  $\mathcal{G}^{\mathcal{R}}_{p, \ell}$ at the subpatch level are
  obtained as images of Cartesian discretizations in parameter space
  (Sections~I-3.2 and~\ref{subsec:overl_dec}.). As mentioned in
  pt.~\ref{time-stepping} below, the solution
  $\textbf{e}^{\mathcal{R}}_{p, \ell}$ is evolved independently on
  each subpatch at each stage of the time-stepping procedure, with
  necessary communication of solution values between subpatches
  enforced via exchange and interpolation
  (Section~\ref{subsec:connectivity}), and with enforcement of
  boundary conditions at all subpatch boundary sections that are
  contained within physical boundaries.
\item\label{viscosity} To ensure solution smoothness and to prevent
  the formation of spurious oscillations in the presence of shocks and
  contact discontinuities, the FC-SDNN method relies on assignment of
  smooth and localized artificial viscosity values. To determine the
  support of the artificial viscosity at each time step, the algorithm
  at first obtains an estimate of the degree of smoothness (denoted by
  $\tau$, see Section~I--5.1) of a certain ``proxy variable'', namely,
  the Mach number
  $\Phi(\textbf{e}) = \norm {\mathbf{u}} \sqrt{\frac{\rho}{\gamma
      p}}$, whose smoothness and/or discontinuity is used as a proxy
  of the corresponding discontinuities in the solution
  $\textbf{e}$. The degree of smoothness $\tau$ of $\Phi$, in turn, is
  estimated on the basis of an underlying, previously trained ANN
  (Section~\ref{sec:viscosity}). Utilizing the degree-of-smoothness
  $\tau$, a preliminary subpatch-wise viscosity assignment I--(38) is
  determined (cf. Section~I--5.1), and scaled, at each point, by an
  estimate of the maximum wave-speed at that point, given by the
  Maximum Wave Speed Bound operator (MWSB) defined in I--(37). To
  ensure spatial smoothness of the artificial viscosity assignments,
  the overall artificial viscosity operator
  $\widehat\mu[\textbf{e}^{\mathcal{R}}_{p, \ell}]$ (equation I--(44))
  incorporates smoothly blended versions of the preliminary
  subpatch-wise viscosity values previously obtained throughout
  $\Omega$. The operator $\widehat\mu$ produces the smoothly-blended
  viscosity values by means of certain smooth windowing functions defined in
  I--(42)---which form a partition of unity within $\Omega$
  subordinated to the covering of $\Omega$ by the set of all
  subpatches---that are evaluated by the algorithm at each point of
  each discretization $\mathcal{G}^{\mathcal{R}}_{p, \ell}$.
\item\label{filtering} In order to control error growth in unresolved
  high frequency modes, at each time step $t_n > 0$ the FC-SDNN method
  applies a spectral filter that multiplies the Fourier coefficients
  by a factor that decays exponentially with the mode number (equation
  I--(32)). The filter is applied independently on each subpatch to
  the $\rho$, $\rho u$ and $\rho v$ components of the solution vector
  $\textbf{e}$ as well as, indirectly, to the energy $E$ (by means of
  the relation between $E$ and the temperature $\theta$), which allows
  a consistent enforcement of adiabatic boundary conditions as
  detailed in Section~I--4.3.1. Additionally, the aforementioned
  exponential filter is used to provide a localized
  discontinuity-smearing strategy for initial data, which is only used
  at the first time step $t_n = 0$---in order to avert the formation
  of spurious oscillations emanating from discontinuous initial
  conditions. This localized initial-data filtering, in turn, relies
  on a somewhat more strongly filtered version of the initial
  condition, which is only used near discontinuities---by smoothtly
  blending, on the basis of a smooth partition of unity, the filtered
  values near discontinuities together with the given unfiltered
  values away from discontinuities; see equation I--(33) and
  Section~I--4.3.2.
\item\label{time-stepping} At each time step $t_n$ the FC-SDNN
  algorithm evaluates the adaptive time step $(\Delta t)_n$ (equation
  I--(31)) and it computes the solution at the next time-step
  $t_{n+1} = t_n + (\Delta t)_n$ by execution of each stage of the
  SSPRK-4 time stepping scheme~\cite{gottlieb2005high} (see
  Section~I--4.1) on the basis of FC-based spatial differentiations
  which incorporate the required Dirichlet and/or Neumann boundary
  conditions (Section~I--4.2) and exchange/interpolation of solution
  values between subpatches
  (Section~\ref{subsec:connectivity}). Patch/subpatch communication is
  imposed at every stage of the SSPRK-4 scheme.
\end{enumerate}

\begin{algorithm}
  \begin{algorithmic}[1]

    \For {every MPI rank $n_r$}
    \State\label{init_grid} Initialize the complete patch/subpatch structure (Section~\ref{subsec:decomposition}), detailing explicit mapping functions, the subpatch decomposition, the subpatch discretization sizes, and the subpatch numbering, but not including assignments of arrays of discretization points. (The resulting data repetition across ranks is not memory-demanding and is thus used.) 
    \State Assign the subpatches $\Omega^{\mathcal{R}}_{p, \ell}$ to the various ranks so as to ensure as close to equidistribution of subpatches per rank as possible.
     \For {every subpatch $\Omega^{\mathcal{R}}_{p, \ell}$ distributed to rank $n_r$}
     \State Allocate the data arrays necessary for storage of discretization points, viscosity values, stage solution values, solution and viscosity data communication.
     \State Compute the physical grids  $\mathcal{G}^{\mathcal{R}}_{p, \ell}$  (Sections~\ref{square} and~\ref{L}) using the mappings on the parameter grid points $(q^1_i, q^2_j) \in G^{\mathcal{R}}_{p, \ell}$ (that are solely determined by the patch/subpatch decomposition and the discretization sizes initialized in line~\ref{init_grid}).
     \State Allocate, compute and store the values of the Jacobian $J_{\mathbf{q}}^{\mathcal{R}, p}$ of the inverse mapping $\big(\mathcal{M}^{\mathcal{R}}_{p}\big)^{-1}$ over all grid points in $\mathcal{G}^{\mathcal{R}}_{p, \ell}$.
     \State Load the trained ANN weights and biases.
     \State  Initialize the unknown solution vector $\textbf{e}_h = \textbf{e}^{\mathcal{R}}_{h, p, \ell}$ (Section~\ref{Preliminaries}) to  the given initial-condition values over all spatial discretization grids $\mathcal{G}^{\mathcal{R}}_{p, \ell}$.
     \State Compute (through communication with the other subpatches) the values of the normalized multi-patch windowing function $\widetilde{\mathcal{W}}^{\mathcal{R}}_{p, \ell, i, j}(\mathbf{x})$  for all $\mathbf{x} \in \mathcal{G}^{\mathcal{R}}_{p, \ell}$.
    
    \EndFor
    \EndFor

  \caption{Multi-patch FC-SDNN initialization. See point~\ref{grids} in Section~\ref{main-elms}}
\label{alg_init}
\end{algorithmic}
\end{algorithm}

\begin{algorithm}
  \begin{algorithmic}[1]
    \For {every MPI rank $n_r$}
        \For {every subpatch $\Omega^{\mathcal{R}}_{p, \ell}$ distributed to rank $n_r$}
             \State Evaluate the proxy variable $\bm{\phi}$ corresponding to $\textbf{e}_h$ at all spatial grid points.
             \State Obtain the smoothness classification values $\tilde \tau[\bm{\phi}]$.
             \State Evaluate the MWSB operator $\widetilde S [\textbf{e}_h]$ at all spatial grid points.
             \State Evaluate the artificial viscosity operator  $\widehat\mu[\textbf{e}^{\mathcal{R}}_{p, \ell}]$.
        \EndFor
    \EndFor
      
  \caption{Parallel evaluation of the preliminary artificial viscosity operator.  See point~\ref{viscosity} in Section~\ref{main-elms}}
\label{alg_prelim_artvisc}
\end{algorithmic}
\end{algorithm}

\begin{algorithm}
  \begin{algorithmic}[1]
    \For {every MPI rank $n_r$}
        \For {every subpatch $\Omega^{\mathcal{R}}_{p, \ell}$ distributed to rank $n_r$}
            \State (Case $t = 0$) Apply localized discontinuity-smearing to the solution vector $\textbf{e}^{\mathcal{R}}_{h, p, \ell}$ and overwrite $\textbf{e}^{\mathcal{R}}_{h, p, \ell}$ with the resulting values.
            \State (Case $t > 0$) Apply the subpatch-wise filtering strategy to the solution vector $\textbf{e}^{\mathcal{R}}_{h, p, \ell}$ and overwrite $\textbf{e}^{\mathcal{R}}_{h, p, \ell}$ with the resulting values. 
        \EndFor
    \EndFor      
  \caption{Parallel filtering of the solution vector $\textbf{e}^{\mathcal{R}}_{h, p, \ell}$.  See point~\ref{filtering} in Section~\ref{main-elms}}
\label{alg_filt}
\end{algorithmic}
\end{algorithm}

\begin{algorithm}
  \begin{algorithmic}[1]
    \For {every MPI rank $n_r$}
    \For {every subpatch $\Omega^{\mathcal{R}}_{p, \ell}$ distributed to rank $n_r$}
    \State Compute the spatial gradient of the viscosity $\mu^{\mathcal{R}}_{p, \ell}$ via FC-based differentiation (Section~\ref{Preliminaries}).
    \State Compute the first- and second-order spatial partial derivative of the density $\rho$, horizontal and vertical momenta $\rho u$ and $\rho v$ and  temperature $T$ via FC-based differentiation while enforcing the relevant Dirichlet or Neumann boundary conditions.
    \State Using the derivatives obtained in lines~3 and~4 together with the Jacobian and the product differentiation  rule, evaluate the discrete differential operator $L^{\mathcal{R}}_{p, \ell}$ by combining these derivatives and the Jacobian $J^{\mathcal{R}, p}_{\textbf{q}}$ of the inverse mapping $\big(\mathcal{M}^{\mathcal{R}}_{p}\big)^{-1}$.
            \State Evaluate the next stage of the SSPRK-4 scheme.
        \EndFor
    \EndFor      
  \caption{Parallel stage evolution and enforcement of boundary conditions for the solution vector $\textbf{e}^{\mathcal{R}}_{h, p, \ell}$. See point~\ref{time-stepping} in Section~\ref{main-elms}}
\label{alg_stage}
\end{algorithmic}
\end{algorithm}

\section{\label{sec:Parallel}Multi-domain parallelization strategy}

This section presents the proposed parallel implementation of the
FC-SDNN algorithm outlined in Section~\ref{main-elms}. The description
assumes the algorithm is run in a number $N_r$ of parallel MPI ranks,
where each rank (i.e., each individual parallel MPI process)
is assumed to be pinned to a single compute core. Of course the code
can be run in serial mode, simply by selecting $N_r =1$.

The proposed computational methodology is encapsulated in
Algorithms~\ref{alg_prelim_artvisc} through~\ref{alg_stage};
Algorithms~\ref{alg_total} and~\ref{alg_init}, in turn, describe the
FC-SDNN driver program, responsible for orchestrating the overall
multipatch MPI implementation of the FC-SDNN solver, and the solver
initialization program, which includes domain initialization as
specified in Item~\ref{grids} of Section~\ref{main-elms}. We note that
the structure of the overall algorithm, which utilizes multiple
patches and subpatches as described in Section~\ref{subsec:overl_dec},
is designed to support arbitrarily large numbers of subpatches within
a fixed patch decomposition of a given geometry. Each subpatch
operates independently of the others, with only minimal communication
required between neighboring subpatches.  This design serves the
triple purpose of (a)~Enabling applicability to general geometries;
(b)~Limiting the size of the required Fourier-continuation expansions;
and, (c)~Facilitating efficient parallelization in a distributed
computing environment. The benefits described in points~(a) and (b)
follow directly by construction. The benefits mentioned in Point (c),
in turn, can be appreciated by examination of
Algorithms~\ref{alg_prelim_artvisc} through~\ref{alg_stage}, all
of which are inherently parallelizable.

To enhance readability we provide the following concise outlines of
the main algorithms, namely, Algorithms~\ref{alg_prelim_artvisc}
through~\ref{alg_stage}.
\begin{enumerate}[(i)]
\item Algorithm~\ref{alg_prelim_artvisc} concerns the evaluation of the subpatch-wise
  preliminary viscosity
  $\widehat\mu[\textbf{e}^{\mathcal{R}}_{p, \ell}]$ which, per
  Item~\ref{viscosity} in Section~\ref{main-elms}, requires the
  computation of (a)~The proxy variable
  $\bm{\phi}^{\mathcal{R}}_{p, \ell}$; (b)~The smoothness
  classification values
  $\tilde\tau[\bm{\phi}^{\mathcal{R}}_{p, \ell}]$; and, (c)~The MWSB
  operator $\widetilde S[\textbf{e}_h]$.
\item Algorithm~\ref{alg_filt} tackles the subpatch-wise construction
  and filtering of all required Fourier Continuation (FC) expansions
  (Item~\ref{filtering} in Section~\ref{main-elms}). This includes:
  (a) Filtering the FC expansions of the solution
  \(\textbf{e}^{\mathcal{R}}_{h, p, \ell}\) (which are computed at
  each time step \(t_n > 0\) and used to produce FC-based spatial
  derivatives, subsequently applied in the SSPRK-4 time-stepping
  procedure~\cite{gottlieb2005high}); and (b) Filtering the FC
  expansions of the initial data to produce localized discontinuity
  smearing, employing specific window functions.
\item Algorithm~\ref{alg_stage} time-steps the solution by executing each stage of
  the SSPRK-4 scheme (Item~\ref{time-stepping} in
  Section~\ref{main-elms}) using FC-based spatial differentiations.
\end{enumerate}
The parallel execution of these algorithms requires communication of
certain solution and viscosity values
(Sections~\ref{subsec:connectivity} and~\ref{sec:viscosity},
respectively) as well as the common value of the variable time-step
$\Delta t$.

\section{\label{sec:numerical results}Numerical results}

This section presents computational results produced by means of the
multi-patch FC-SDNN method in a number of challenging test cases,
including problems involving interactions between
supersonic/hypersonic flow and obstacles, and including multiple
moving-shocks and contact discontinuities, as well as shock collisions
with obstacles, domain boundaries, contacts, and other shocks,
etc. The efficiency of the parallel implementation of the algorithm,
in terms of both weak and strong scaling, is demonstrated in
Section~\ref{subsec:scaling}, and a number of illustrative numerical
examples produced by the multi-patch FC-SDNN algorithm are presented
in Section~\ref{subsec:applications}.

The numerical tests presented in this chapter were conducted on a
Beowulf computer cluster comprising 30 dual-socket nodes, each
equipped with 384 GB of GDDR4 RAM. All 30 nodes are interconnected via
HDR Infiniband. Each socket contains two 28-core Intel Xeon Platinum
8273 processors, for a total of 56 cores per node and 1680 cores over
the 30-node cluster. While supported by the Xeon processors, the
hyper-threading capability was not utilized in any of the numerical
examples presented in this section. The ``production run'' shown in
Figure~\ref{schlieren_M10_Shockmatrix} required approximately 472 CPU
hours, while those in Figure~\ref{Prism_Shock_solutions} took
approximately 125 CPU hours each.  Further details on runtimes and
parallel efficiency are presented in Section~\ref{subsec:scaling}.

\subsection{\label{subsec:scaling}Parallel performance: Weak and
  strong scaling}


This section demonstrates the weak and strong parallel scalability
enjoyed by our implementation of the multi-patch FC-SDNN algorithm. In
particular, Sections~\ref{scalingIII} and~\ref{scalingIV} present weak
scaling results: in Section~\ref{scalingIII} the size and complexity
of the physical problems considered are increased, whereas in
Section~\ref{scalingIV} the discretizations are increased (refined)
for fixed physical problems. Section~\ref{scalingV}, in turn,
demonstrates the strong scaling of the algorithm employing the physical
problem considered in Section~\ref{scalingIV}. The two types of physical
problems considered are described in Section~\ref{scalingI}, and the
parallel scalability metrics used are introduced in
Section~\ref{scalingII}.

\subsubsection{Scaling I: Test problems\label{scalingI}}
\begin{figure}[H]
  \centering
  \includegraphics[width=1\linewidth,]{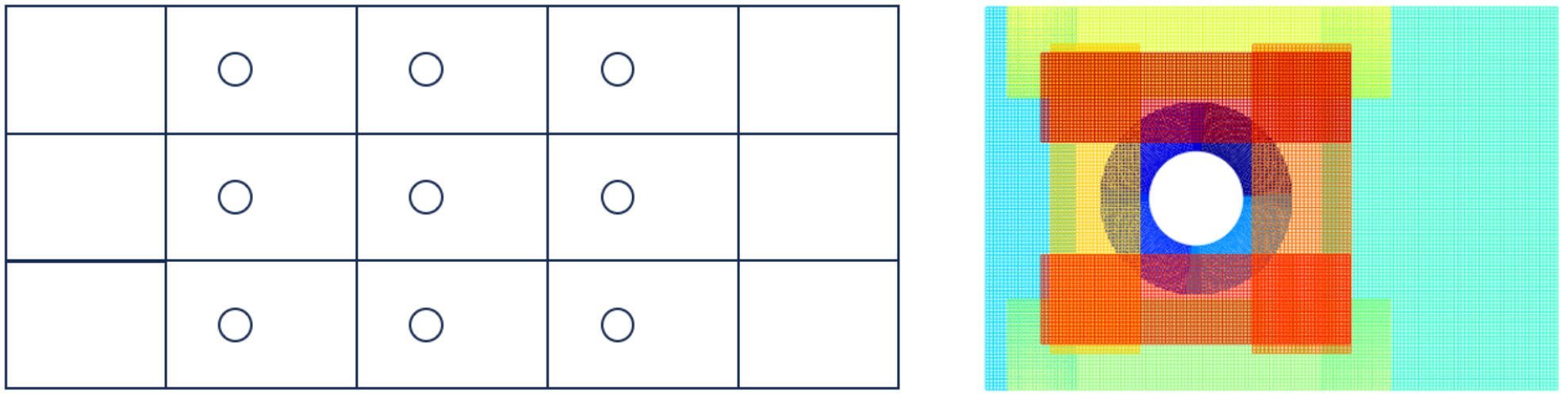}
  \caption{Illustration of the types of geometries used for the test
    problems in Section~\ref{scalingI}. Left panel: initial partition
    of the computational domain utilized in Test Problem~1 with
    $n_\mathrm{row} = 3$ and $n_\mathrm{col} = 3$. Right panel: patch
    decomposition of each cylinder-containing subdomain in the left
    panel.  The right panel also depicts the complete computational
    domain used for Test problem~2.}
  \label{fig:MatrixCylinderDomain}
\end{figure}

\paragraph{Test Problem 1: Mach 10 shock mitigation by matrices of
  cylindrical obstacles.} Test Problem 1 actually comprises a set of
problems in which a Mach 10 shock-wave (a shock wave traveling at 10
times the speed of sound of the unperturbed fluid) impinges upon a
rectangular array of circular cylinders containing $n_\mathrm{col}$
columns and $n_\mathrm{row}$ rows of cylinders, as depicted in the
left panel of Figure~\ref{fig:MatrixCylinderDomain}. The computational
domain considered is divided into three zones, namely, an
obstacle-free front region ahead of the rectangular array; a middle
region containing the rectangular array of obstacles; and an
obstacle-free wake region. The middle
$n_\mathrm{col}\times n_\mathrm{row}$-cylinder region is constructed
as a corresponding array of identical overlapping rectangular
subdomains of horizontal and vertical sides of lengths $\ell_x = 3$
and $\ell_y = 2$, respectively, each one of which consists of $54$
subpatches arranged around a single cylinder of radius $0.25$. One
such rectangular subdomain is depicted in the right panel of
Figure~\ref{fig:MatrixCylinderDomain}. The front and wake regions
comprise a vertical array of $n_\mathrm{row}$ rows of rectangular
cylinder-free patches of horizontal and vertical sides of lengths
$L_x = 2.5$ and $L_y = 2$, respectively, each one of which once again
contains $54$ rectangular subpatches.

The problem prescriptions are completed by incorporating the following
initial and boundary conditions. The initial conditions utilized are
given by
\begin{equation} \label{IC_Mach10_shock} (\rho, u, v, p)= \left\lbrace
    \begin{array}{ccc}
      (4.0816, 8.25, 0, 116.5) & \mbox{if} & x \leq 0.6 \\
      (1.4, 0, 0, 1) & \mbox{if} & x > 0.6. \\          
    \end{array}\right.
\end{equation}
An inflow condition with $(\rho, u, v, p)$ values coinciding with the
$x\leq 0.6$ initial values, on the other hand, is imposed at the left
boundary at all times. An outflow condition consisting of the
time-independent pressure value $p = 1$ is imposed at the right
boundary.  Slip-wall (zero-normal velocity) boundary conditions are
imposed at the bottom and top walls. No slip (zero velocity) and
adiabatic (zero normal component of the gradient of the temperature)
boundary conditions, finally, are imposed at the boundaries of the
cylinders at all times---since, as detailed as part of the
Test-Problem-2 description below, a viscous problem is solved near the
cylinder boundaries whenever non-zero velocities occur at the cylinder
boundaries.

\paragraph{Test Problem 2: Hypersonic flow past a cylindrical
  obstacle.} Test Problem 2 concerns a hypersonic (Mach 10) flow past
a cylinder of diameter $0.25$ and centered at $(x_c, y_c) = [1.1, 0]$,
where the computational domain corresponds exactly to a single
instance of a subdomain containing a cylinder as previously described
as part of Test Problem 1 and shown in the right panel of
Figure~\ref{fig:MatrixCylinderDomain}. An initial configuration with
$54$ subpatches is considered, which is then refined in accordance
with the method described in Section~\ref{subsubsec:subpatches} in
order to perform weak and strong scaling tests. The initial conditions
assumed correspond to a Mach 10 flow given by
\begin{equation}\label{IC_Mach10_flow}
  (\rho, u, v, p) = (1.4, 10, 0 ,1).
\end{equation}
As in Test Problem~1, an inflow condition with $(\rho, u, v, p)$
values coinciding quantitatively with the initial values is imposed at
the left boundary at all times, and no boundary conditions are imposed
on the right outflow boundary, as befits a supersonic
outflow. Reflecting boundary conditions, corresponding to zero-normal
velocity, are imposed at the bottom ($y = -1$) and top ($y = 1$)
walls. No slip (zero velocity) and adiabatic (zero normal component of
the gradient of the temperature) boundary conditions are imposed at
the boundaries of the cylinder at all times---as befits the
viscous-like problem that is solved in a neighborhood of cylinder on
account of the artificial viscosity which, as illustrated in
Figure~6 in Part~I, is assigned by the
artificial-viscosity algorithm in that region.

\subsubsection{Scaling II: Parallel performance metrics\label{scalingII}}
In order to quantify the parallel scaling efficiency of the FC-SDNN
algorithm in multi-patch settings, for given numbers $N$ of
discretization points and $N_C$ of compute cores, and using a number
$N^0_C$ of compute cores for reference, we denote by $T_S(N_C, N)$ and
$T_S(N^0_C, N)$ the number of seconds required by the FC-SDNN to
advance $4 N$ unknowns for one time-step using $N_C$ cores and $N^0_C$
cores, respectively. The strong (resp. weak) scaling efficiencies
$E^s_{N^0_C, N_C}$ (resp. $E^w_{N^0_C, N_C}$) are then defined by
\begin{equation}\label{scaling}
E^s_{N^0_C, N_C} = \frac{T_S(N^0_C, N) N^0_C}{T_S(N_C, N) N_C}, \qquad E^w_{N^0_C, N_C} = \frac{T_S(N^0_C, N)} {T_S(N_C, N N_C / N^0_C)}.
\end{equation}
An alternative measure of the parallel computing time required by the
algorithm to advance the $4N$ unknowns associated with an $N$-point
discretization grid in a given computational experiment is provided by
the number $S_{N_C}$ of CPU-seconds required to advance the simulation
for one time-step in a number $N_C$ of CPU cores per $10^6$ unknowns,
that is,
\begin{equation}\label{Sm}
S_{N_C, N} = \frac{N_C \times \mbox{(total computation time)} \times 10^6}{4 \times N  \times \mbox{(Time steps)}}.
\end{equation}
In particular, the strong and weak efficiencies can be re-expressed as
a function of $S_{N_C}$:
\begin{equation}\label{scaling_sm}
E^s_{N^0_C, N_C} = \frac{S_{N^0_C, N} N^0_C}{S_{N_C, N} N_C}, \qquad E^w_{N^0_C, N_C} = \frac{S_{N^0_C, N}} {S_{N_C, N_c / N^0_C N}}.
\end{equation}

Section~\ref{scalingIII} utilizes Test Problem 1, described in
Section~\ref{scalingI}, to demonstrate the algorithm's weak scaling in
scenarios where the number of discretization points is increased to
tackle progressively larger physical problems. In this context, the
problem sizes are expanded by incrementally increasing the numbers of
rows and/or columns of cylindrical obstacles in the rectangular array
of cylinders considered.  Section~\ref{scalingIV} then utilizes Test
Problem 2 in Section~\ref{scalingI} to illustrate once again weak
scaling properties of the algorithm, but this time in a context in
which increasing discretizations result from mesh refinement. Test
Problem 2 is used a second time in Section~\ref{scalingV}, finally, to
study strong scalability---by solving a fixed physical problem defined
on a fixed mesh on the basis of increasingly larger numbers of
computing cores.

The very high, essentially perfect, weak scalability demonstrated by
the tests presented in Sections~\ref{scalingIII} and~\ref{scalingIV}
indicate that both the communication and the average computing-cost
per subpatch required by the algorithm remain essentially fixed as the
numbers of cores and the sizes of the problems increase
proportionally. The high but less-than-perfect strong scalability
illustrated in Section~\ref{scalingV}, on the other hand, reflects the
fact that, owing to differences in communication costs, the computing
time per subpatch may vary among subpatches. The strong scaling could
be further improved by incorporating an algorithm that distributes a
combination of subpatch types to each computer node, including
adequate proportions of high and low communication-cost
subpatches. Such additional algorithmic developments are beyond the
scope of this paper, however, and are left for future work.

\subsubsection{\label{scalingIII}Scaling III: Weak scaling under
  problem enlargement}


This section illustrates the weak scalability of the FC-SDNN algorithm
by enlarging the problem size via a progressive addition of columns
and/or rows of cylindrical obstacles in a rectangular array of cylinders
of the type described in Test Problem 1, Section~\ref{scalingI}, and
advancing the solution up to time T = 0.1 (with space and time units
such that the speed of sound in the unperturbed flow state is $a = 1$,
and the distance between the centers of two vertically consecutive
cylinders is 1.25).  As showcased in Table~\ref{table:weak2}, the
FC-SDNN algorithm enjoys essentially perfect weak scaling as the
number of discretization points and cores are proportionally increased
($N/N_C = N_\mathcal{Q}$ with $N_\mathcal{Q} = 101^2$, see
equation~\eqref{eq:NQ} and associated text), starting with an
$N_C = 9\times 54 = 486$ computing-core initial configuration (at 54
cores per node) for a Test Problem-1 array containing a single
three-cylinder column in addition to the front and wake regions.

\begin{table}[H]
\centering
\begin{tabular}{|c|c|c|c|c|c|c|c|}
\hline
$N$ & $n_\mathrm{row}$ & $n_\mathrm{col}$ & $N_C$ & $\mbox{Nodes}$ & $T(s)$ &  $S_{N_C}(s)$ & $E^w_{486, N_C} (\%)$ \\ \hline
4,957,686 & 3 & 1 & 486       & 9      &   74.4   & 0.97 &  \\
6,610,248 & 3 & 2  & 648       & 12      &  74.4  & 0.97 & 100 \\
6,610,248 & 4 & 1 & 648       & 12       &  75.0  & 0.98 & 99 \\
8,262,810 & 3 & 3 & 810        & 15       &  74.2  & 0.97 & 100 \\
  8,262,810 & 5 & 1 & 810        & 15       &  74.4  & 0.97 & 100 \\
  8,813,664 & 4 & 2 & 864        & 16       &  74.8  & 0.98 & 99 \\
  9,915,372 & 3 & 4 & 972        & 18       &  74.0  & 0.97 & 101 \\
  11,017,080 & 4 & 3 & 1080        & 20       &  74.0  & 0.97 & 101 \\
  11,017,080 & 5 & 2 & 1080        & 20       &  74.8  & 0.98 & 99 \\
  11,567,934 & 3 & 5 & 1134        & 21       &  74.0  & 0.97 & 101 \\
  13,220,496 & 4 & 4 & 1296        & 24       &  73.9  & 0.97 & 101 \\
  13,771,350 & 5 & 3 & 1350        & 25       &  74.5  & 0.97 & 100 \\
  15,423,912 & 4 & 5 & 1512        & 28       &  74.8  & 0.98 &  99 \\
  16,525,620 & 5 & 4 & 1620        & 30       &  75.0  & 0.98 & 99 \\  \hline  
\end{tabular}
\caption{FC-SDNN weak parallel scaling for various rectangular arrays
  of cylinders (Test Problem~1 in Section~\ref{scalingI}), containing
  $n_\mathrm{row}$ rows and $n_\mathrm{col}$ columns of cylindrical
  obstacles, and using $N_C$ cores, with $N_C$ ranging from $486$ to
  $1620$ ($9$ to $30$ computer nodes).}
\label{table:weak2}
\end{table}

\subsubsection{\label{scalingIV}Scaling IV: Weak scaling under
  mesh refinement}

This section once again illustrates the weak scalability of the
proposed FC-SDNN parallel implementation, but this time in the context
of mesh refinement---wherein, as described in
Section~\ref{subsubsec:subpatches}, increases in the numbers of
discretization points are obtained via corresponding increases in the
number of subpatches used while simultaneously and proportionally
increasing the number $N_C$ of cores---in a setting where the number
of cores equals the number of subpatches. We do this here by solving
Test Problem 2 up to time $T = 0.1$ (with space and time units such
that the speed of sound in the unperturbed flow state is $a = 1$, and
the distance between the top and bottom boundaries of the
computational domain equals $2$). As the mesh is refined the time-step
decreases with the average spatial meshsize in an approximately
proportional fashion (following equation~I--(31)), thereby increasing
in a (roughly) linear manner the number of simulation time-steps---as
illustrated in Table~\ref{table:weak1}.

Table~\ref{table:weak1} displays the runtimes $T$, the numbers
$S_{N_C}$ and the efficiencies $E^w_{54, N_C}$ as the total number
of cores and discretization points are increased proportionally, as
noted in the table, using 54 cores per computer node and increasing
the number of discretization points $N$ proportionally to the number
of computer nodes used.

\begin{table}[H]
\centering
\begin{tabular}{|c|c|c|c|c|c|c|}
\hline
$N$ & $N_C$ & $\mbox{Nodes}$ & $T(s)$     & $\mbox{Time steps}$ & $S_{N_C}(s)$ & $E^w_{54, N_C} (\%)$ \\ \hline
550,854 & 54     & 1 & 23.0       & 563      &   1.00      &   \\
2,203,416    & 216   & 4  & 44.41       & 1106      &  0.98   & 102 \\
4,957,686    & 486    & 9 & 64.96       & 1650       &  0.96   & 104 \\
8,813,664    & 864     & 16 & 86.20        & 2194       &  0.96  & 104 \\
13,771,350    & 1350    & 25 & 107.46        & 2738       &  0.96  & 104 \\ \hline
\end{tabular}
\caption{FC-SDNN weak parallel scaling in a mesh refinement context
  (Test Problem 2 in Section~\ref{scalingI}) using $N_C$ cores, with
  $N_C$ ranging from $54$ to $1350$ ($1$ to $25$ computer nodes).}
\label{table:weak1}
\end{table}

As illustrated in Table~\ref{table:weak1}, the weak scaling efficiency
of the procedure in the present mesh-refinement context is excellent,
steadily rising above $100\%$ with respect to the coarsest mesh. We
attribute this better than perfect scalability to slightly diminishing
workloads associated to the ranks carrying the heaviest workloads as
the mesh is refined. In detail, for initial, coarse, subpatch
partitions, a number of subpatches along the boundary of a patch play
the roles of both donors and recipients of interpolated solution
data. Since the patches (and their overlaps) remain unchanged as the
subpatch refinements take place, donors tend not to be receivers and
vice versa, and thus the maximum interpolation-data workload per MPI
rank that occurs in such patches is decreased. This process eventually
stabilizes (at $104\%$ of the efficiency for the initial $N$ value),
as illustrated in Table~\ref{table:weak1}, once most subpatches
engaged in communication of interpolated data only play the roles of
either donor or recipient.


\begin{table}[H]
\centering
\begin{tabular}{|c|c|c|c|c|c|c|}
\hline
$N$ & $N_C$ & $\mbox{Nodes}$ & $T(s)$  & $S_{N_C}(s)$    & $E^s_{54, N_C} (\%)$ & $E^s_{N_C / 2, N_C} (\%)$     \\ \hline
 & 54   & 1  & 1126.0     & 0.79  &       &         \\
 & 108 & 2   & 598.8    & 0.84   & 94      &  94   \\
8,813,664 & 216  & 4  & 330.1    & 0.92    & 85       &  91   \\
 & 432  & 8   & 167.0   & 0.93    & 84       &  99   \\
 & 864  & 16  & 86.4    & 0.96    & 81       &  97  \\ \hline
\end{tabular}
\caption{FC-SDNN strong parallel scaling (Test Problem 2 in
  Section~\ref{scalingI}) using $N_C$ cores, with $N_C$ ranging from
  $54$ to $864$ ($1$ to $16$ computer nodes).}
\label{table:strong}
\end{table}

\subsubsection{\label{scalingV}Scaling V: Strong scaling}

This section illustrates the strong parallel scalability of the
FC-SDNN algorithm in a context in which the number of allocated
processing cores is increased for a given physical problem and for
fixed patch/subpatch decomposition of the problem---namely Test
Problem 2 discretized on the basis of $864$ subpatches, which
corresponds to a number $N = 8,813,664$ of discretization
points. (This example thus utilizes the same physical problem
considered in Section~\ref{scalingIV}, using the next-to-finest
discretization considered in that section, but assigning multiple
subpatches to each computer core, as needed for each given number of
computers cores used to deal with the fixed number of subpatches
considered.)  Strong scaling is studied in the 30-node cluster used by
progressively increasing the number of nodes utilized to solve the
problem up to time $T=0.1$, and using a constant number of $54$ cores
per node. The associated times $T(s)$ in seconds and parallel
efficiencies $E^s_{54, N_C}$ with respect to a run of $54$ cores are
displayed in Table~\ref{table:strong}. The observed decrease in the
strong scaling efficiency as the number of cores is increased may be
attributed to differences in communication costs and associated
computational workloads assigned to various subpatches, as discussed
towards the end of Section~\ref{scalingII}.
\subsection{\label{subsec:applications}Applications}

This section presents numerical results of the application of the
FC-SDNN to a number of physical problems of various degrees of
geometric complexity, including results concerning strong shocks as
well as supersonic and hypersonic flows. The new results clearly
accord with previous numerical simulations, theory and experimental
data, and they include simulations at high Mach numbers, including
e.g. Mach 10 shocks and flows---for which good agreement is observed
with flow details predicted by previous theory as well as
extrapolations from experimental data such as e.g. the distance
between an obstacle and a reflected bow shock.

As indicated above, a variety of numerical results are presented in
this section. Thus, to assess the impact of the multi-patch
decomposition used, Section~\ref{validation} compares multi-patch
results to results obtained by the single-patch version of the FC-SDNN
method, and Section~\ref{engy_cons} presents results concerning energy
conservation in the multi-patch
setting. Section~\ref{subsubsec:initialization} then presents a
description of the initial conditions concerning problems considered
in the subsequent sections, namely, supersonic flows past an unbounded
triangular wedge (Section~\ref{subsubsec:wedge}) and a bounded
triangular prism (Section~\ref{subsubsec:flow_prism}), shock-wave
mitigation by arrays of obstacles (Section~\ref{Matrixcylinders}), and,
finally, shock past a prism (Section~\ref{subsubsec:shock_prism}).

To facilitate visualization of shock waves,
following~\cite{banks2008study} and~\cite{nazarov2017numerical},
Schlieren diagrams, namely, displays of the quantity
\begin{equation}\label{exp_schlieren}
  \sigma = \exp \left( -\beta \frac{\lvert \nabla \rho(\mathbf{x}) \rvert - \min_{\mathbf{x} \in \Omega} \lvert \nabla \rho(\mathbf{x}) \rvert }{\max_{\mathbf{x} \in \Omega}\lvert \nabla \rho(\mathbf{x}) \rvert - \min_{\mathbf{x} \in \Omega} \lvert \nabla \rho(\mathbf{x}) \rvert} \right)
\end{equation}
with $\beta = 10$, are presented for each test case. The time-stepping
constant $\textrm{CFL} = 0.25$ in~I--(31) was used for all the
simulations involving $C_1$-type corners and patches considered in
Sections~\ref{subsubsec:flow_prism} and~\ref{subsubsec:shock_prism}
(since these problems involve rather stretched meshes and thus small
distances between spatial discretization points), while the value
$\textrm{CFL} = 0.5$ was used in all other cases.

\subsubsection{\label{validation} Multi-patch method validation}
In order to validate the accuracy of the multi-patch FC-SDNN algorithm
we consider the Riemann4 test case introduced
in~\cite{lax1998solution} and tackled
in~\cite[Sec. 4.3.2]{bruno2022fc} by means of a single-patch FC-SDNN
solver. The left panel in Figure~\ref{fig: R4_multipatch} reproduces
the single-patch $600 \times 600$-point solution previously presented
in the middle panel of Figure~16 in~\cite{bruno2022fc}. The right
panel in Figure~\ref{fig: R4_multipatch}, in turn, displays the
solution produced, for the same problem, on the basis of a
$7 \times 7$ subpatch decomposition with a similar mesh spacing:
$h \approx 0.002$. The two images in Figure~\ref{fig: R4_multipatch}
present closely matching flow features, with errors consistent with
the underlying error level---which may be gleaned from the red error
curves in Figure~11 in~\cite{bruno2022fc}. In particular, the
smoothness of the contour levels is maintained, and the shocks are
equally well resolved in both cases.
\begin{figure}
  \centering
  \includegraphics[width=0.35\linewidth]{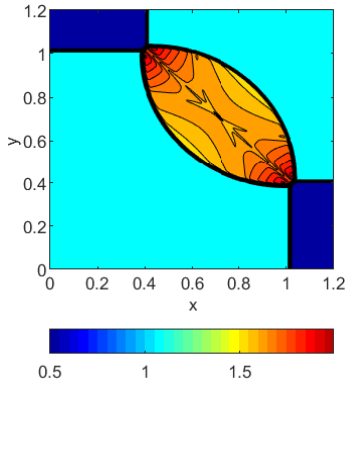}
  \includegraphics[width=0.35\linewidth]{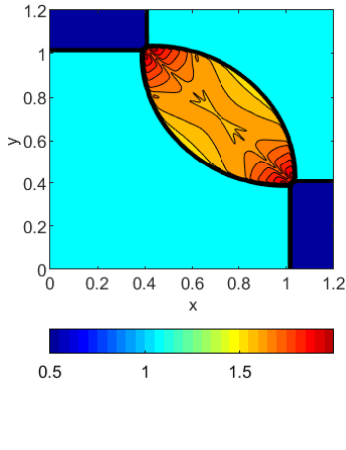}
 \vspace{-1.3cm}
 \caption{FC-SDNN numerical solutions for the Riemann4 test case at
   $t = 0.25$ visualized using thirty equispaced contours between
   $\rho = 0.5$ and $\rho = 1.99$. Left: Single-domain
   simulation~\cite[Fig. 16(b)]{bruno2022fc} using $360,000$
   discretization points. Right: Multi-patch FC-SDNN simulation using
   a single patch subdivided into $7 \times 7$ subpatches, totaling
   approximately $500,000$ discretization points. In both cases, the
   mesh spacing $h$ is $\approx 0.002$.\looseness -1 \label{fig:
     R4_multipatch}}
\end{figure}

\subsubsection{\label{engy_cons} Energy conservation.} Several
mechanisms used by the FC SDNN algorithm, such as the introduction of
artificial viscosity (Section~\ref{sec:viscosity}) to avoid spurious
oscillations, and the use of a spectral filtering at every time-step
to control the error growth in unresolved high frequency modes
(point~3 in Section~\ref{main-elms}), are dissipative in nature, but,
as demonstrated below, do not lead to excessive energy loss. To
quantify the effect we studied a 2D multi-patch version of the Sod
Problem~\cite[Sec. 4.3.1]{bruno2022fc}, in the domain
$[0, 1] \times [0, 0.25]$, with initial conditions given by
\begin{equation}
(\rho, u, v, p)=
\left\lbrace
    \begin{array}{ccc}
        (1, 0, 0, 1) & \mbox{if} & x < 0.5   \\ 
        (0.125, 0, 0, 0.1) & \mbox{if} & x \geq 0.5,
    \end{array}\right.
\end{equation}
with reflecting boundary conditions $u = 0$ at $x = 0$ and $x = 1$,
and with vanishing normal derivatives for all variables on the top and
bottom boundaries $y=0.25$ and $y=0$ at all times---so as to simulate
a vertically infinite domain. As no energy enters or exits the domain
via inflow or outflow, the exact total energy value
$\bar{E}_\mathrm{exact} = \int_0^{0.25}\int_0^1
E_\mathrm{exact}(x,y,t)dxdy$ contained in the domain is constant---as
it can be checked easily e.g. by differentiation with respect to time
of the integral of $E$ followed by use of the energy
equation~\eqref{eq: euler 2d equation} together with vanishing values
of the normal derivatives of $u$ and $v$ at the upper and lower
boundaries $y=0.25$ and $y=0$. It is easy to check that for the exact
solution we have $\bar{E}_\mathrm{exact} = 1.375$ in this case. Thus
any observed departures $\Delta\bar{E}$ of the computed energy
$\bar{E}$ from this exact value provide an indication of the
dissipative character of the algorithm.

\begin{figure}
        \includegraphics[width=0.5\linewidth]{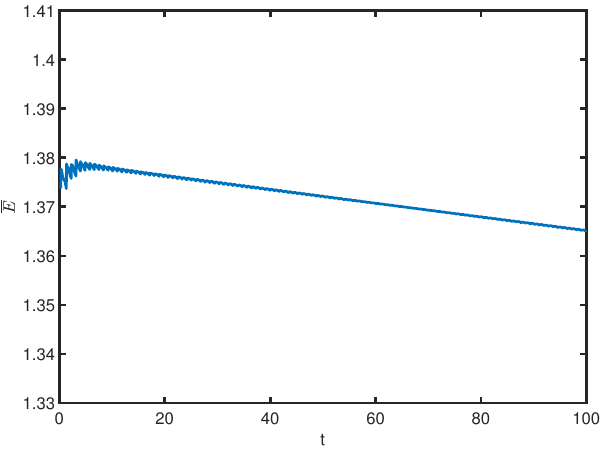}
        \includegraphics[width=0.5\linewidth]{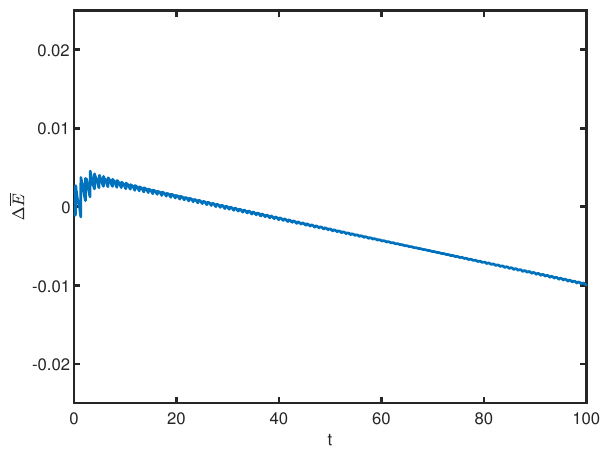}
        \caption{Numerical values of the total energy $\bar{E}$ (left
          panel) and its defect $\Delta\bar{E}$ from the exact value
          $\bar{E}_\mathrm{exact} = 1.375$ (right panel) as functions
          of time in the time interval $0\leq t\leq 100$---within
          which the shock, contact discontinuity and rarefaction wave
          were reflected 60 times per boundary, for a total of 120
          wall reflections. The figures show that the numerical energy
          $\bar{E}$ remains essentially constant
          ($\bar{E} \approx 1.37\pm 0.01$) for the duration of the
          time-interval considered, with defect values consistent with
          the algorithm's error levels illustrated in~\cite[Fig. 11]{bruno2022fc}.}
        \label{fig: SodRefl}
\end{figure}

For the purpose of this energy conservation test, the computational
domain was discretized by using a unique patch decomposed into four
overlapping horizontal subpatches. The solution was computed up to
time $T = 100$, during which the shock was reflected 60 times per
boundary, for a total of 120 wall reflections. As showcased in
Figure~\ref{fig: SodRefl}, the numerical energy values $\bar{E}$
remain essentially constant, $\bar{E} \approx 1.37\pm 0.01$, where the
errors are consistent with the underlying error level indicated by the
red error curves in~\cite[Fig. 11]{bruno2022fc}. Noticeable features
of the energy and energy error curves include an early-time increase
in the energy values (with a total increase that diminishes in size as
the discretization is refined, as it was found on the basis of an
additional set of tests not included in this paper for brevity) as
well as a slow energy decrease attributable to the dissipative
processes under consideration.
    
\subsubsection{\label{subsubsec:initialization}Initial conditions used
  in Sections~\ref{subsubsec:wedge}
  through~\ref{subsubsec:shock_prism}}
Two different kinds of initial conditions are used in what
follows. The first type concerns supersonic/hypersonic ``Mach $M$''
flows past obstacles (for $M$ values such as $1.4$, $3.5$, $10$ etc.),
which are considered in Sections~\ref{subsubsec:wedge}
and~\ref{subsubsec:flow_prism}. The initial conditions
\begin{equation}\label{mach_flow_ic}
(\rho, u, v, p) = (1.4, M, 0 ,1),
\end{equation}
we use in all such cases indeed correspond to a uniform Mach $M$ flow
with speed of sound $a = 1$ (per the identity
$a^2 = \frac{\gamma p}{\rho}$ together with the assumed ideal
diatomic-gas constant $\gamma = 1.4$).  Initial conditions of a second
kind are utilized in Sections~\ref{Matrixcylinders}
and~\ref{subsubsec:shock_prism}, which correspond to the interaction
of a shock initially located at $x = x_s$ and traveling toward the
region $x \geq x_s$ at a speed $M$ that is supersonic/hypersonic with
respect to the speed of sound $a = 1$ ahead of the shock. For such a
Mach-$M$ shock the initial conditions are given by
\begin{equation} \label{mach_shock_ic}
(\rho, u, v, p)=
\left\lbrace
    \begin{array}{ccc}
        (\frac{(\gamma + 1) M^2}{(\gamma - 1) M^2 + 2}, \frac{\zeta - 1}{\gamma M}, 0, \zeta) & \mbox{if} & x < x_s   \\ 
        (1.4, 0, 0, 1) & \mbox{if} & x \geq x_s,
    \end{array}\right.
\end{equation}
where $\zeta = \frac{2 \gamma M^2 - \gamma + 1}{\gamma +1}$ denotes
the strength of the shock; see e.g.~\cite{chaudhuri2011use}.

\subsubsection{\label{subsubsec:wedge}Supersonic flow past a triangular wedge}

\begin{figure}
  \begin{subfigure}[t]{0.49\linewidth}
    \centering
    \includegraphics[width=1\linewidth]{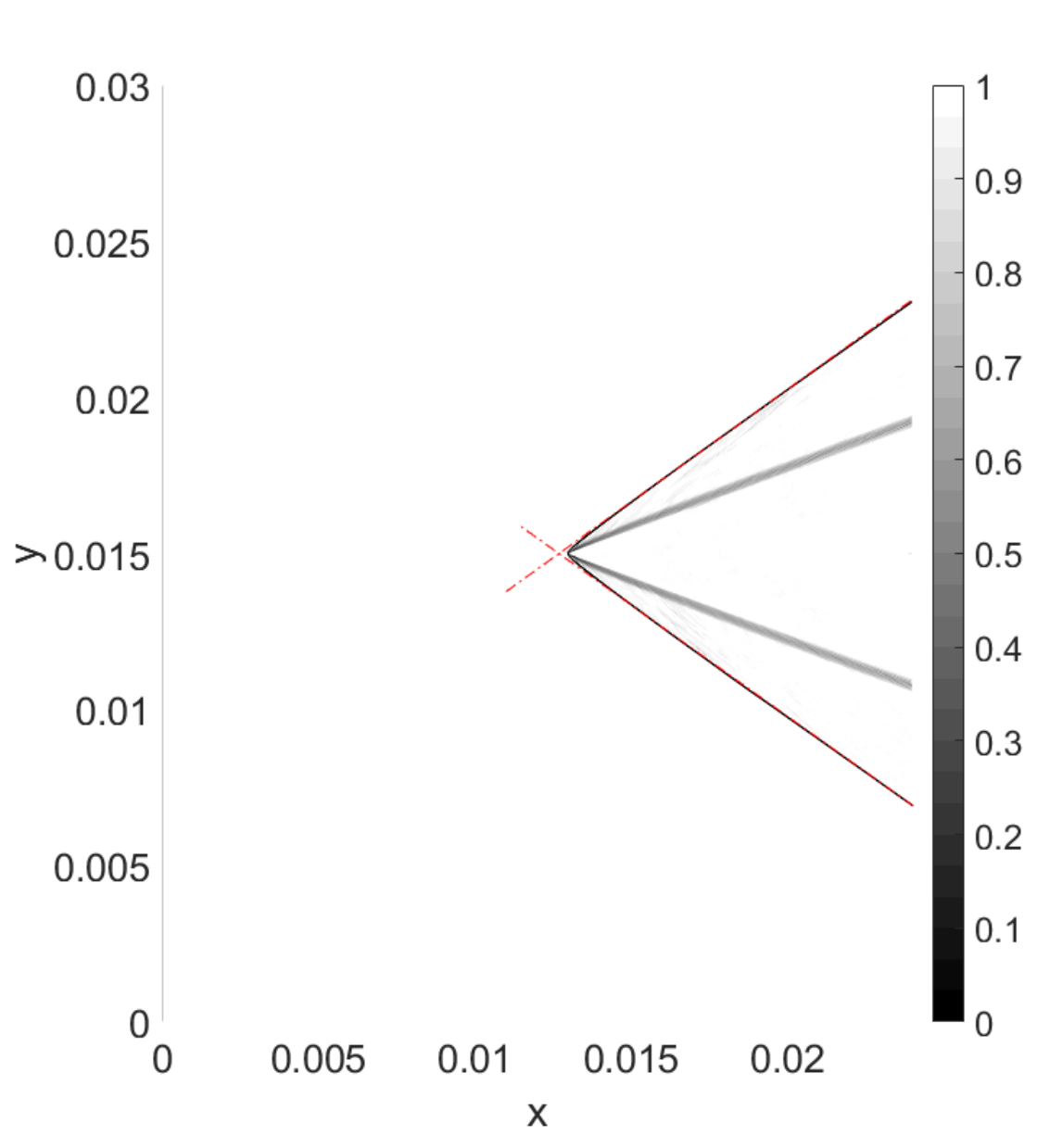}
  \end{subfigure} 
  \begin{subfigure}[t]{0.49\linewidth}
    \centering  
    \includegraphics[width=1\linewidth]{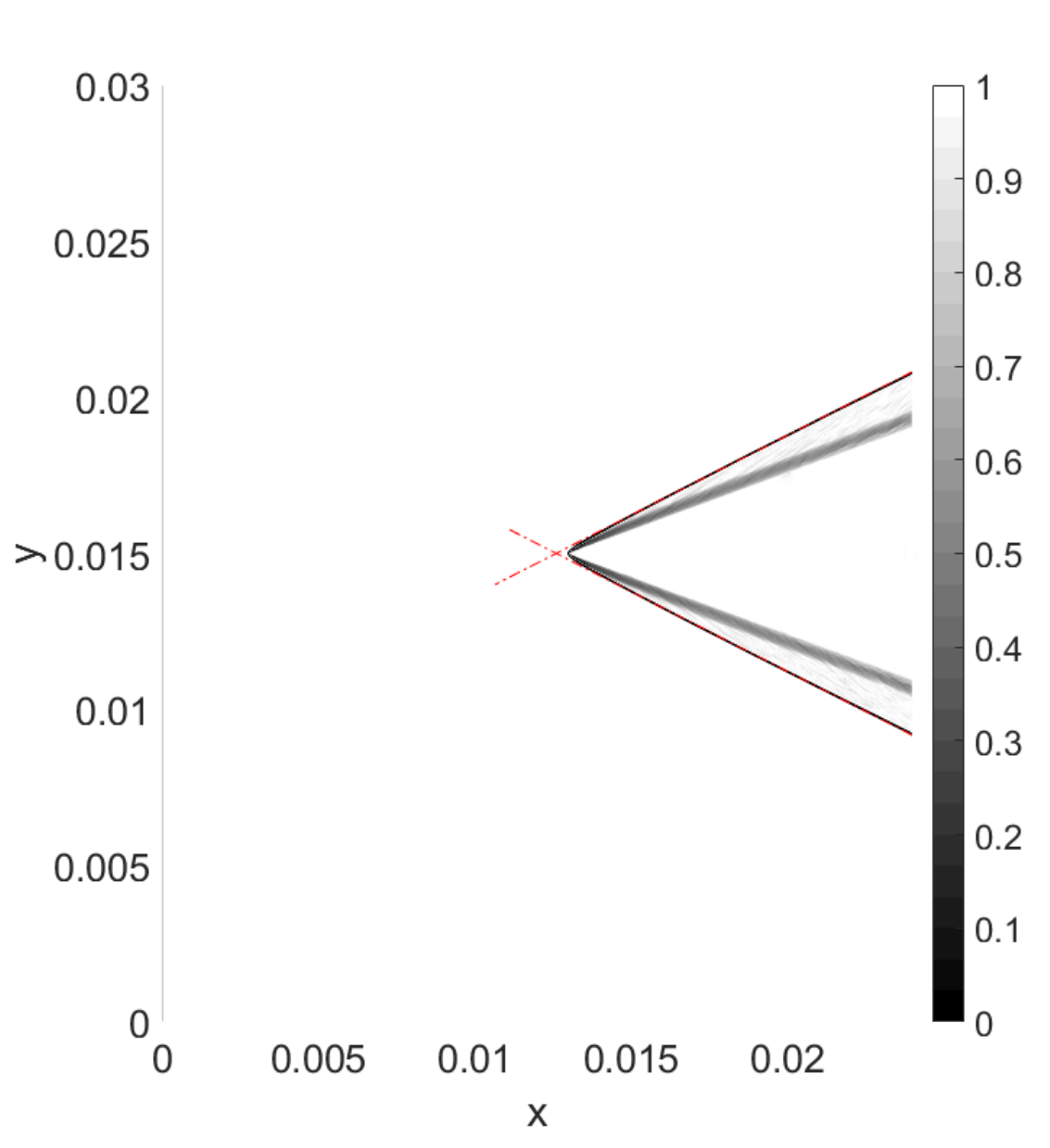}
  \end{subfigure}

  \caption{Supersonic flows past a triangular wedge at time
    $T = 0.02$, shown for two different Mach numbers:, Mach 3.5 (left)
    and Mach 10 (right). Solutions were computed on a mesh with
    $N = 31.5M$ discretization points. In both cases, red lines
    indicate linear fits to the shock fronts.}
\label{Wedge_Flow_solutions}
\end{figure}

We next consider supersonic and hypersonic flows past the triangular
wedge depicted in Figure~\ref{Wedge_Flow_solutions}, whose tip is
located at the point $(x_t, y_t) = (0.013, 0.015)$, and whose interior
angle equals $\alpha_t = 40^{\circ}$. The computational domain
consists of the portion of the rectangle $[0, 0.024] \times [0, 0.03]$
located outside the wedge. Initial conditions given by
equation~\eqref{mach_flow_ic} with supersonic Mach number $M = 3.5$
and hypersonic Mach number $M = 10$ are considered in what follows. An
inflow condition with $(\rho, u, v, p)$ values coinciding
quantitatively with the initial values is imposed at the left boundary
at all times, and no boundary conditions are imposed on the right
boundary, as befits a supersonic outflow. Vanishing normal derivatives
for all variables are imposed at the top on bottom domain
boundaries. No slip and adiabatic boundary conditions are imposed on
the wedges at all times. The results at time $T = 0.02$ for Mach
$M = 3.5$ and $M = 10$ flows and for a mesh containing approximately
31 million points are presented in
Figure~\ref{Wedge_Flow_solutions}.

Notably, a straight oblique shock forms starting at the tip of the
triangular wedge, with a deflection angle $\alpha_d$ with respect to
the horizontal. A closed form relation between $\alpha_d$, the wedge
angle $\alpha_t$ and the Mach number $M$ of the incoming flow can be
obtained in closed form (for an infinite wedge) on the basis of the
Rankine-Hugoniot jump conditions~\cite{chaudhuri2011use}; the result
is
\begin{equation}\label{deflection}
\tan(\frac{\alpha_t}{2}) = 2 \cot(\alpha_d) \frac{M^2 \sin^2(\alpha_d) - 1}{M^2(\gamma + \cos(2 \alpha_d)) + 2}.
\end{equation}
The numerical results presented in Figure~\ref{Wedge_Flow_solutions}
show a deflected shock which curves around the tip of the wedge to
eventually form two straight shocks, one on each side of the
obstacle. The near-tip shock curvature can be attributed to the
formation of a thin viscous boundary layer along the wedge boundaries,
consistent with the use of artificial viscosity in~\eqref{viscosity}
and the accompanying no-slip and adiabatic boundary conditions imposed
on the obstacle, as described earlier in this section. The deflection
angle corresponding to the computational simulation was calculated by
fitting a straight line (red dashed line in
Figure~\ref{Wedge_Flow_solutions}) to the shock starting from
$x = 0.024$ (the rightmost abscissa shown in the figure) which, in
view of the aforementioned slight near-tip curvature, meets the $x$
axis at an abscissa $x_d$ slightly ahead of the tip abscissa $x_t$.
Table~\ref{table:deflection} presents a comparison between the
numerically computed deflection angles for various meshes and the
corresponding theoretical values. As demonstrated in that table, the
deflection angle  approaches the expected theoretical value as
the discretization is refined. (The use of uniformly refined meshes,
which is required by the present implementation of our algorithm,
limits the maximum resolution achievable in near-tip mesh within a
reasonable computing time, but it is expected that use of temporal
subcycling and/or implicit near-tip solvers would enable significantly
higher resolution, and thus, an even closer approximation of the
theoretical value, within reasonable computational times.)  The
$M=3.5$ shock problem was previously considered
in~\cite{chaudhuri2011use}; that contribution calculates a deflection
angle of $34.5^\circ$ on the basis of a mesh of approximately one
million discretization points, which more closely approximates the
theoretical deflection angle than any of the predictions presented in
Table~\ref{table:deflection}. We suggest that further work is
warranted in this connection however, as the coarser discretization
used in~\cite[Fig. 4]{chaudhuri2011use} results in a wider shock
profile and an associated imperfect match to a straight line showing
noticeable departures from the shock for extended shock sections near
the wedge vertex.

\begin{table}
\centering
\begin{tabular}{|c|c|c|}
\hline
  & $M = 3.5$ &  $M = 10$ \\ \hline
$\mbox{Inviscid theory}$ & $34.6^{\circ}$ &  $25.8^{\circ}$  \\
$N = 520,251$ &  $36.8^{\circ}$   & $21.9^{\circ}$\\
  $N = 2,356,431$ & $36.0^{\circ}$   & $23.0^{\circ}$  \\
  $N = 10,007,181$ & $35.7^{\circ}$ &  $23.8^{\circ}$ \\
  $N = 31,500,688$ & $35.2^{\circ}$   & $24.3^{\circ}$ \\ \hline
\end{tabular}
\caption{Deflected shock angle for the problem of a supersonic flow
  past a wedge, including inviscid-theory values and computational
  results for various numbers $N$ of discretization points.}
\label{table:deflection}
\end{table}

\subsubsection{\label{subsubsec:flow_prism}Supersonic flow past a triangular prism}

\begin{figure}[H]
    \centering
    \includegraphics[width=0.4\linewidth]{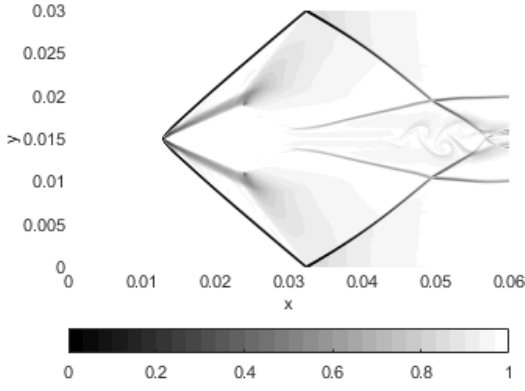}\hspace{1cm}
    \includegraphics[width=0.4\linewidth]{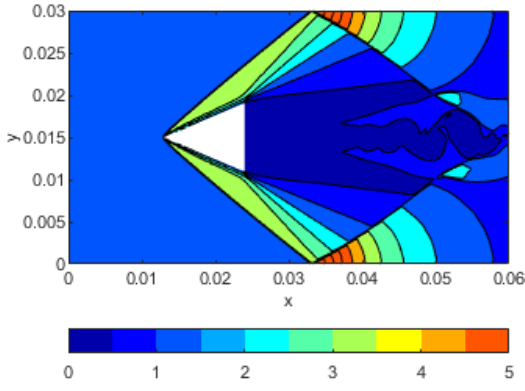}
    \includegraphics[width=0.4\linewidth]{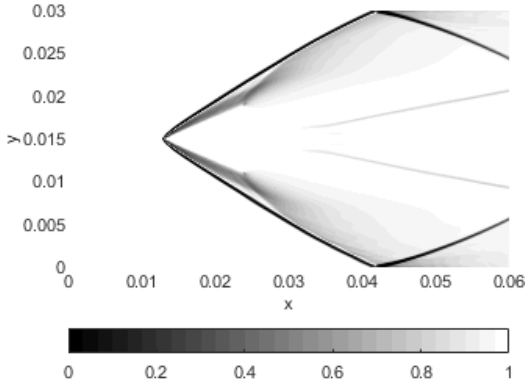}\hspace{1cm}
    \includegraphics[width=0.4\linewidth]{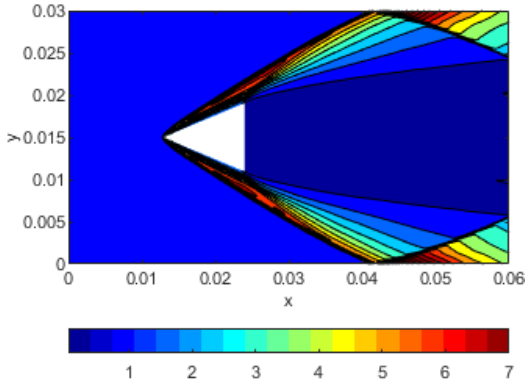}
    \caption{Supersonic flow past a triangular prism at time
      $T = 0.02$, shown for Mach 3.5 (top row) and Mach 10 (bottom
      row). Both simulations were performed on a mesh containing
      $N = 10.5$ million points. Left column: density Schlieren
      images. Right column: density contour plots.}
\label{Prism_Flow_solutions}
\end{figure}

This section concerns flow problems in a wind tunnel similar to the
one considered in~\cite{chaudhuri2011use}, of dimensions
$[0, 0.06]\times[0, 0.03]$, over a stationary triangular prism of
length $\ell = 0.011$ from left to right, whose front vertex has
coordinates $(x_t, y_t) = (0.013, 0.015)$ and front angle
$\alpha_t = 40^{\circ}$. Two different initial value problems are
considered for this geometry, namely those resulting from supersonic-
and hypersonic-flow initial conditions given by
equation~\eqref{mach_flow_ic} with Mach numbers $M = 3.5$ and
$M = 10$, respectively. Reflecting boundary conditions are imposed at
the bottom and top walls, together with no slip and adiabatic boundary
conditions at the prism's boundaries. An inflow condition with
$(\rho, u, v, p)$ values coinciding quantitatively with the initial
values is imposed at the left boundary at all times, and no boundary
conditions are imposed on the right boundary, as befits a supersonic
outflow. Schlieren images as well as a contour plots of the density
$\rho$ for both test cases at time $T = 0.02$ are presented in
Figure~\ref{Prism_Flow_solutions}. As in
Section~\ref{subsubsec:wedge}, oblique shocks are reflected off the
front-facing sides of the triangular prism. As shown in the figures,
these shocks are reflected by the top and bottom wind-tunnel
boundaries, and the reflected shocks, in turn, interact with the wake
at the back of the prism. Clearly, the FC-SDNN method allows for a
fine resolution of the shock structures, and as demonstrated by the
contour plots in Figures~\ref{Prism_Flow_solutions}, smooth density
profiles are preserved away from shocks. This problem presents certain
challenges concerning the possible formation of vacuum states in the
wake of the prism, in particular for high Mach number flows, but,
unlike previous spectral, finite-element and finite-difference-based
approaches, no density limiters were used to preserve density
positivity. To the best of the authors' knowledge, the numerical
simulations of Mach 3.5 flows past a prism reported
in~\cite{chaudhuri2011use,daspreliminary} are the fastest wind-tunnel
flow-past-prism numerical simulations presented to date.

\subsubsection{\label{Matrixcylinders}Shock-wave mitigation by a matrix of solid cylindrical obstacles}

\begin{figure}
    \centering
\includegraphics[width=1\linewidth,]{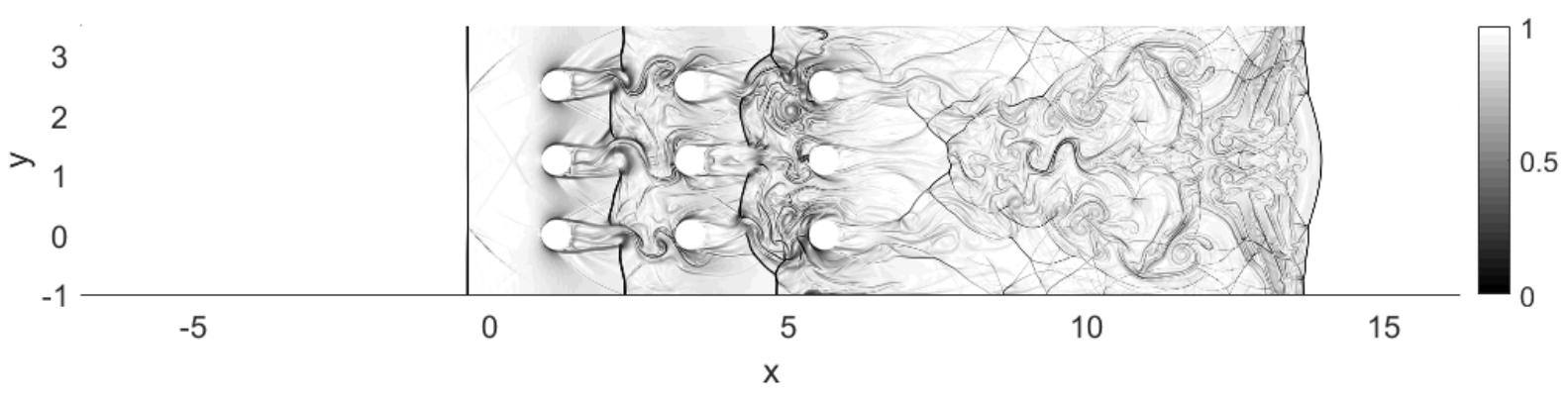}
\caption{Density Schlieren visualization of a Mach 10 shock-wave
  mitigation by a 3x3 matrix of cylindrical obstacles at time
  $T = 1.4$, obtained on an $N \approx 16.5M$-point mesh.}
    \label{schlieren_M10_Shockmatrix}
  \end{figure}

This section presents results for a Mach 10 shock propagating through
an array of cylindrical obstacles---a significantly stronger shock
than those previously studied for this type of configuration: prior
work concerns Mach 1.8 and Mach 3 shocks in similar
geometries~\cite{chaudhuri2013numerical,mehta2016numerical}. The study
of such shock-mitigation configurations is important for the design of
effective shielding systems~\cite{chaudhuri2013numerical}. The
geometry considered, depicted in
Figure~\ref{fig:MatrixCylinderDomain}, results as a slight
modification of the one used in Test Problem 1
(Section~\ref{scalingI}) with $n_\mathrm{row} = 3$ and
$n_\mathrm{col} = 3$; the modifications introduced here concern the
front and the wake region, both of which have been extended: the new
front (resp. wake) region considered here contains a $3 \times 3$
array of rectangular cylinder-free patches (resp. a $4 \times 3$ array
of cylinder-free patches), instead of the corresponding $1\times 3$
front and wake arrays considered in Section~\ref{scalingI}.  The
initial conditions considered here are given by
equation~\eqref{mach_shock_ic} with shock traveling at Mach numbers
$M = 3$ and $M = 10$, and starting at the position $x_s = 0.5$; the
boundary conditions used are identical to those described in Test
Problem 1. Figure~\ref{schlieren_M10_Shockmatrix} presents a Schlieren
diagram of the solution at time $T = 1.4$.

\subsubsection{\label{subsubsec:shock_prism}Shock past a prism}

\begin{figure}[H]
    \centering
    \includegraphics[width=0.4\linewidth]{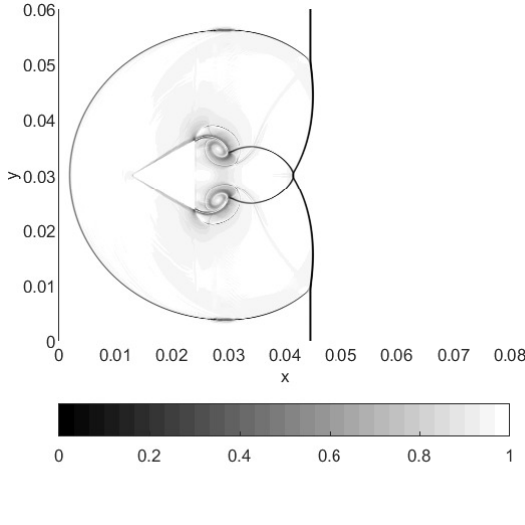}\hspace{1cm}
    \includegraphics[width=0.4\linewidth]{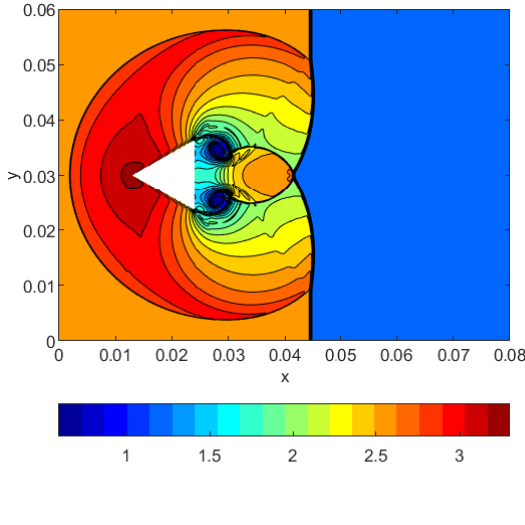}
    \includegraphics[width=0.4\linewidth]{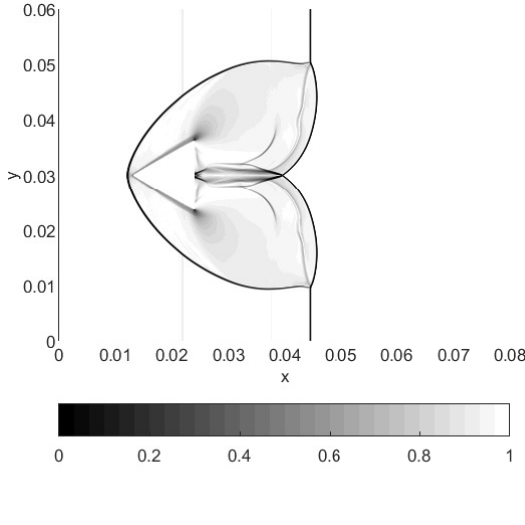}\hspace{1cm}
    \includegraphics[width=0.4\linewidth]{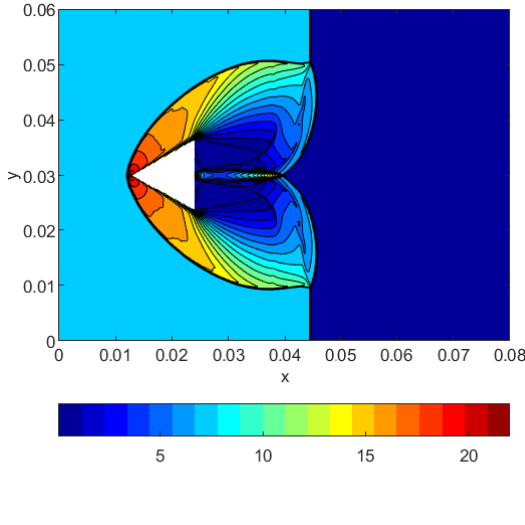}\vspace{-0.7cm}
    \caption{Shock–prism interaction at two different shock speeds and
      times: Mach 1.5 at $T = 0.02$ (top) and Mach 10 at $T = 0.00375$
      (bottom). Both simulations were performed on a mesh with
      approximately $N = 10.5$ million points. Left column: density
      Schlieren images. Right column: density contour plots.}
\label{Prism_Shock_solutions}
\end{figure}

We finally consider the interaction of a shock with a stationary
triangular prism in a $[0, 0.08]\times[0, 0.06]$ wind tunnel. In this
case an equilateral triangular prism of length $\ell = 0.011$ from
left to right is considered, with a front vertex at coordinates
$(x_t, y_t) = (0.013, 0.03)$ and front angle
$\alpha_t = \frac{\pi}{3}$. The initial conditions for the two tests
considered here are given by equation~\eqref{mach_shock_ic} with
$M = 1.5$ and $M = 10$, respectively, and with $x_s = 0.007$. An
experimental test for a similar geometry, the \textit{Schardin}
problem, was introduced in~\cite{schardin1957high} for a Mach 1.34
shock, and studied numerically in~\cite{chang2000shock} and
in~\cite{chaudhuri2011use}.  An inflow condition with
$(\rho, u, v, p)$ values coinciding with the $x\leq x_s$ initial
values are imposed at the left boundary at all times. Further, an
outflow condition consisting of the time-independent pressure value
$p = 1$ is imposed at the right boundary.  Slip-wall boundary
conditions are imposed at the bottom and top walls, while no-slip and
adiabatic boundary conditions are imposed at the boundaries of the
prism at all times. Solutions at times $T = 0.025$ for the Mach 1.5
shock and $T = 0.00375$ for the Mach 10 shock are displayed in
Figure~\ref{Prism_Shock_solutions}.  The numerical solution correctly
displays transmitted and reflected shocks, as well as slip lines and
vortices at the back of the prism. As in the previous examples, the
density contour plots in Figure~\ref{Prism_Shock_solutions} exhibit
smooth contour levels away from flow discontinuities.  As in previous
examples, the density contour plots in
Figure~\ref{Prism_Shock_solutions} display smooth contour levels away
from flow discontinuities. Also noteworthy are the contact
discontinuities which, originating at the upper and lower triple
points and, in the $M = 1.5$ cases, extending toward the upper and
lower vortices, appear as clearly visible traces in the Schlieren
plots and as blips in the contour plots. The clear presence of these
structures underscores the low dispersion and dissipation inherent in
the discretization strategy employed, as such features would be
significantly smeared or entirely lost in presence of a higher
numerical dissipation or dispersion.
 
\section{\label{sec:conclusion}Conclusions}

This two-paper sequence introduced a novel multi-patch computational
algorithm for general gas-dynamics problems: the FC-SDNN {\em
  spectral} shock-dynamics solver. The method applies to {\em general
  domains}, it can handle {\em supersonic} and {\em hypersonic} flows,
and, unlike other spectral schemes, it operates without restrictive
CFL constraints. Part~I extended the single-patch FC-SDNN
solver~\cite{bruno2022fc} to general smooth domains and successfully
tackled challenging high-Mach-number flows. Building on this
foundation, Part II focused on: 1)~Handling domains with corners
through the introduction of $\mathcal{C}_1$ and $\mathcal{C}_2$-type
corner patches; 2)~Introducing the algorithm's parallel
implementation, demonstrating, in particular excellent strong parallel
scaling, and near-perfect weak scaling; 3)~Validating the multi-patch
method by comparing its solutions with those produced by the
single-patch solver on simple domains and by demonstrating its energy
conservation properties; and, 4)~Benchmarking the method against
theoretical and experimental results for complex geometries with
corners. In particular this two-paper sequence has demonstrated that
the proposed multi-patch strategy delivers for general domains the
benefits previously achieved in~\cite{bruno2022fc} on the basis of
single spectral expansions on simple domains.  In view of its smooth
but localized viscosity assignments this procedure effectively
eliminates Gibbs oscillations while avoiding introduction of the
flow-field roughness that is often evidenced by the serrated level
sets produced by other methods. In view of its use of FC-based Fourier
expansions, further, the proposed algorithm exhibits spectral behavior
away from shocks (resulting, in particular, in limited dispersion in
such regions~\cite[Sec. 4.1.2]{bruno2022fc}) while enabling the
evaluation of solutions in general domains and under general boundary
conditions. Unlike other techniques, finally, the approach does not
require use of problem-dependent algorithmic parameters.

A variety of extensions are envisioned as part of
future work, including the adaptation of the method to general 3D
configurations, potentially incorporating CAD-described
(Computer-Aided Design) structures. Other possible developments, which
would be valuable for improving efficiency and accuracy, include the
design of hybrid ADI-based implicit/explicit versions of the
scheme~\cite{bruno2019higher} and/or the use of temporal subcycling
(e.g. to achieve high resolution near edges and corners without
requiring uniformly refined meshes), and the efficient
GPU parallelization to enhance computational performance.


\section*{Acknowledgments}
The authors gratefully acknowledge support from NSF and AFOSR under
contracts DMS-2109831, FA9550-21-1-0373 and FA9550-25-1-0015.  \small

\small

\bibliographystyle{abbrv}
\bibliography{bibliography}

\end{document}